\documentclass[12pt]{article}
\usepackage{fullpage,graphicx,psfrag,amsmath,verbatim,xcolor}
\usepackage[small,bf]{caption}

\usepackage[noend]{algpseudocode}
\usepackage{algorithm}
\usepackage{url}
\usepackage{bbm}
\usepackage{tikz}
\usepackage{tikz}
\usetikzlibrary{positioning,calc,shapes,arrows}
\tikzstyle{block} = [rectangle, draw,
text width=7em, text centered, minimum height=5em]
\tikzstyle{line} = [draw, -latex',line width=1mm]

\newcommand{\ones}{\mathbf 1}

\newcommand{\reals}{{\mbox{\bf R}}}

\newcommand{\BEAS}{\begin{eqnarray*}}
\newcommand{\EEAS}{\end{eqnarray*}}
\newcommand{\BEA}{\begin{eqnarray}}
\newcommand{\EEA}{\end{eqnarray}}
\newcommand{\BEQ}{\begin{equation}}
\newcommand{\EEQ}{\end{equation}}
\newcommand{\BIT}{\begin{itemize}}
\newcommand{\EIT}{\end{itemize}}

\newcommand{\diag}{\mathop{\bf diag}}

\newcommand{\geomean}{\mathop{\bf geomean}}

\newcommand{\argmin}{\mathop{\rm argmin}}

\newcommand{\dom}{\mathop{\bf dom}} 

\newcommand{\sign}{\mathop{\bf sign}}

\newcommand{\eg}{{\it e.g.}}
\newcommand{\ie}{{\it i.e.}}

\newcounter{algorithmctr}[section]
\renewcommand{\thealgorithmctr}{\thesection.\arabic{algorithmctr}}
\newenvironment{algdesc}%
   {\refstepcounter{algorithmctr}
   \begin{list}{}{%
       \setlength{\rightmargin}{0\linewidth}%
       \setlength{\leftmargin}{.05\linewidth}}%
       \rmfamily\small
       \item[]{\setlength{\parskip}{0ex}\hrulefill\par%
        \nopagebreak{\bfseries\textsf{Algorithm \thealgorithmctr~}}}}%
   {{\setlength{\parskip}{-1ex}\nopagebreak\par\hrulefill} 
   \end{list}}

   {\refstepcounter{algorithmctr}
   \begin{list}{}{%
       \setlength{\rightmargin}{0\linewidth}%
       \setlength{\leftmargin}{.05\linewidth}}%
       \rmfamily\small
       \item[]{\setlength{\parskip}{0ex}\hrulefill\par%
        \nopagebreak{\bfseries\textsf{Algorithm}}}}%
   {{\setlength{\parskip}{-1ex}\nopagebreak\par\hrulefill} 
   \end{list}}

\title{Implementation of an\\ Oracle-Structured Bundle 
Method\\ for Distributed Optimization}

\author{Tetiana Parshakova \and Fangzhao Zhang 
\and Stephen Boyd}

\begin{document}
\maketitle

\begin{abstract}
We consider the problem of minimizing a function that is a sum of convex agent
functions plus a convex common public function that couples them. The
agent functions can only be accessed via a subgradient oracle; the public function is
assumed to be structured and expressible in a domain specific language (DSL) for
convex optimization.  We focus on the case when the evaluation of the agent
oracles can require significant effort, which justifies the use of solution methods
that carry out significant computation in each iteration.
To solve this problem we integrate multiple known techniques 
(or adaptations of known techniques)
for bundle-type algorithms, obtaining a method
which has a number of practical advantages over other methods that are compatible
with our access methods, such as proximal subgradient methods.
First, it is reliable, and works well across a number of applications.  Second,
it has very few parameters that need to be tuned, and works
well with sensible default values.
Third, it typically produces a reasonable approximate solution in just a few tens
of iterations.
This paper is accompanied by an open-source implementation of the
proposed solver, available at \url{https://github.com/cvxgrp/OSBDO}.
\end{abstract}

\section{Oracle-structured distributed optimization} \label{s-prob}

\subsection{Oracle-structured optimization problem}
We consider the optimization problem 
\BEQ\label{e-prob}
\begin{array}{ll}
\mbox{minimize} & h(x) = f(x) + g(x),
\end{array}
\EEQ
with variable $x \in \reals^n$, where $f,g:\reals^n \to \reals \cup \{ \infty\}$
are the oracle and structured objective functions, respectively.
We assume the problem has an optimal point $x^\star$, and
we denote the optimal value of the problem \eqref{e-prob} as $h^\star=h(x^\star)$.
We use infinite values of $f$ and $g$ to 
encode constraints or the domains, with 
$\dom f = \{ x\in \reals^n \mid f(x)<\infty\}$ and similarly for $g$.
We assume that $\dom f \supseteq \dom g \neq \emptyset$,
\ie, every point in the domain of $g$ is in the domain of $f$,
and $g$ has at least one point in its domain.

\paragraph{The oracle objective function.}
We assume that the oracle function $f$ is block separable,
\[
f(x) = \sum_{i=1}^M f_i(x_i), \quad x=(x_1, \ldots, x_M),
\]
with $f_i:\reals^{n_i}\to \reals \cup \{\infty\}$ closed convex,
where $n_1 + \cdots + n_M = n$.
We refer to $x_i$ as the variable and $f_i$ as the objective function of agent $i$.
Our access to $f_i$ is only via an oracle that evaluates $f_i(x_i)$ and
a subgradient $q_i \in \partial f_i(x_i)$ at any point $x_i \in \dom f_i$.

\paragraph{The structured objective function.}
We assume the structured function $g$ is closed convex.
While $f$ is block separable, $g$ (presumably) couples the 
block variables $x_1, \ldots, x_M$.
We assume that $g$ is given in a structured form, using the language of
Nesterov, which means we have a complete description of it.  We assume
that we can minimize $g$ plus some additional structured function of $x$.
As a practical matter this might mean that $g$ is expressed in a 
domain-specific language (DSL) for convex optimization, such as 
CVXPY~\cite{diamond2016cvxpy,agrawal2018rewriting},
based on disciplined convex programming (DCP) \cite{grant2006disciplined}.

\paragraph{Example.}
Problem \eqref{e-prob} is very general, and includes as special cases
many convex optimization problems arising in applications.
To give a simple and specific example, it includes the so-called
consensus problem \cite[Ch.~7]{boyd2011distributed},
\[
\begin{array}{ll}
\mbox{minimize} & \sum_{i=1}^M f_i(x_i)  \\
\mbox{subject to} & x_1 = \cdots = x_M,
\end{array}
\]
where $x_i \in \reals^p$ are the variables, and $f_i$ are given convex functions.
We put this in the form \eqref{e-prob} with $n_i=p$, $n=Mp$, $x=(x_1, \ldots, x_M)$,
and 
\[
g(x) = \left\{ \begin{array}{ll} 0 &
x_1 = \cdots = x_M\\
\infty & \mbox{otherwise},
\end{array} \right.
\]
the indicator function of the consensus constraint $x_1 = \cdots = x_M$.

\paragraph{Optimality condition.}
The optimality condition for problem \eqref{e-prob} is: $x$ is optimal if
and only if
\BEQ\label{e-opt-cond} 
\partial h(x) = \partial f(x) + \partial g(x) \ni 0, 
\EEQ
where $\partial h(x)$ denotes the subdifferential of $h$ at $x$.
In some applications we may be interested in finding a subgradient
$q^\star\in \partial f(x^\star)$ for which $-q^\star \in \partial g(x^\star)$.
Such a subgradient can sometimes be interpreted as a vector of optimal prices.
The method we describe in this paper will also compute (an estimate of) $q^\star$.

\paragraph{Our focus.}
We seek an algorithm that solves the distributed oracle-structured problem 
\eqref{e-prob}, respecting our access assumptions.  
Several generic methods can be used, such as 
proximal subgradient methods or their accelerated extensions,
described below.
Most research on methods for the composite minimization problem
focus on the case where $f$ is differentiable and $g$ has a simple,
typically analytically computable, proximal operator,
and algorithms that involve very minimal computation beyond evaluation of 
the gradient of $f$ and the proximal operator of $g$, such as a 
few vector-vector operations.
Our focus here, however, is on the case where
the agent oracles can be expensive to evaluate, which has several 
implications.  
First, it means that we focus on algorithms that 
in practice find good suboptimal points in relatively few iterations.
Second, it justifies algorithms that solve an optimization problem
involving $g$ in each iteration, instead of carrying out 
just a few vector-vector operations.

\paragraph{Our contribution.}
Our contribution is to assemble a number of known methods, such as 
diagonal preconditioning, level bundle methods, and others into
an algorithm that works well on a variety of practical problems,
with no parameter tuning.  
By working well, we mean that modest 
accuracy, say on the order of $1\%$, is achieved typically in tens 
of iterations.
In the language of the bundle method literature, our final algorithm
is a disaggregate partially exact bundle method.

Our open-source implementation, along with all data to reproduce 
all of our numerical experiments, is available at 
\begin{quote}
\url{https://github.com/cvxgrp/OSBDO}.
\end{quote}
OSBDO stands for oracle-structured bundle distributed optimization.

\subsection{Previous and related work}

There is a vast literature on general distributed optimization, but 
fewer authors consider the specific subgradient oracle plus structured function
access we consider here.
Several general methods can be used to solve the problem \eqref{e-prob} using our
access methods, including 
subgradient, cutting-plane, and bundle-type methods.

\paragraph{Subgradient methods.}
Subgradient methods were originally developed by
Shor and others in the 1970s \cite{shor2012minimization}. 
Early work that describes subgradients and convex optimization 
includes \cite{rockafellar1981TheTO}. (To use a subgradient-type
method for our problem, we would need to 
compute a subgradient of $g$, which is readily done.)  
Subgradient methods typically 
require a large number of iterations, and are employed when the 
computational cost
of each iteration is low, which is not the case in our setting.
Subgradient methods also involve many algorithm parameters that must be
tuned, such as a step-size sequence. 
Modern variations include AdaGrad \cite{duchi2011adaptive}, an adaptive subgradient 
method which shows good practical results in the online learning setting,
but for our setting still requires far too many iterations to achieve even
modest accuracy.

\paragraph{Proximal subgradient methods.}
A closely related method that is better matched to our specific access 
restrictions is the proximal subgradient 
method, which in each iteration requires a subgradient of $f$
and an evaluation of the proximal operator of $g$.
(This proximal step is readily carried out since $g$ is structured.)
Most of the original work on these types of methods focuses on the case
where $f$ is differentiable, and the method is called the proximal gradient
method.  The proximal gradient method can be described as an 
operator-splitting method~\cite[Ch.~4.2]{parikh2014proximal}; some early work includes
\cite{bruck1975iterative, chen1997convergence}.
Since then there have been two relevant developments.
The papers \cite{passty1979ergodic, lions1979Splitting} handle the case where $f$
is nondifferentiable, and the gradient is replaced with a subgradient,
so the method is called the proximal subgradient method.
The stochastic case is addressed in \cite{schechtman2022stochastic},
where a stochastic proximal subgradient method is developed.
Another advance is the development of simple
generic acceleration methods, originally 
introduced by Nesterov \cite{nesterov1983amf}.
Other work addresses issues such as
inexact computation of the proximal operator \cite{birgin2003inexact},
or inexact computation of the subgradient \cite{burachik2015anas}.
Proximal gradient and subgradient methods
are widely used in different applications; see, \eg, the book
\cite{combettes2011proximal}, which covers proximal
gradient method applications in signal processing.

Compared to the method we propose, proximal subgradient methods fail to
take into account previously
evaluated function values and subgradients (other than through their effect on
the current
iterate), and our ability to build up a model of each agent function separately.
These are done in the method we propose, which requires only a modest 
increase in complexity over evaluating
just the proximal operator of $g$, and give substantial improvement in 
practical convergence.

\paragraph{Cutting-plane methods.} Cutting-plane methods
can be traced back to Cheney and Goldstein
\cite{cheney1959newton} and Kelley's cutting-plane method \cite{kelley1960cutting}.
These methods maintain a piecewise affine lower bound or minorant on
the objective, and improve it in each iteration using the subgradient and value
of the function at the current iterate.  Each iteration requires the solution of a 
linear program (LP), with size that increases with iterations. 
Cutting-plane methods are extended to handle convex mixed-integer problems in
\cite{westerlund1995extended}.
Limited-memory or constraint-dropping versions, that drop terms
in the minorant so the LP solved in each iteration does not grow in size,
are given in 
\cite{elzinga1975central, dem1985nondifferentiable}. 
Constraint dropping for general outer approximation algorithms are also 
considered in \cite{gonzaga1979constraint}.
Many other variations on cutting-plane
methods have been developed, including the analytic center 
cutting-plane method introduced by
Atkinson and Vaidya \cite{atkinson1995acp}. 
A review of cutting-plane methods used in machine learning can be found in
\cite{sra2012optimization}.

\paragraph{Bundle methods.}
Bundle methods are closely related to cutting-plane methods.
The first difference is the addition of a stabilization term, among 
which a proximal regularization term is most common and leads to 
the so-called proximal bundle method~\cite{kiwiel1990proximity, frangioni2020standard}. 
Typically in each iteration of the proximal bundle 
method a quadratic program (QP) must be solved.  Other stabilization forms 
include the level bundle method \cite{lemarechal1995new, frangioni2020standard} 
and the trust-region bundle method \cite{10.2307/169692, frangioni2020standard}.
The second difference between cutting-plane and bundle methods 
is the logic that updates the current point only if 
a certain sufficient descent condition holds.

Bundle methods were first developed as dual 
methods in \cite{lemarechal1975aneo} and \cite{mifflin1977semismooth}; 
the primal form of bundle methods was mainly studied in the 1990s. 
The first convergence proof for the proximal bundle method is given in
\cite{lemarechal1978nonsmooth}.
A comprehensive review of the history and development of bundle methods 
can be found in \cite[Ch.~XIV,XV]{hiriart1996convex}. 
Later variations of bundle methods include variable
metric bundle methods which accumulate second order information about 
function curvature in the proximal regularization 
term \cite{mifflin1996quasi, chen1999proximal, luksan1999globally, 
burke2000ots, haarala2004new}, bundle methods that handle inexact
function values or subgradient values 
\cite{hintermuller2001proximal, kiwiel1985algorithm, kiwiel1995approximations, 
kiwiel2006proximal, hare2016proximal, van2017probabilistic, lv2018proximal,
oliveira2020bundle}, 
bundle methods with semidefinite cutting sets replacing the traditional 
polyhedral cutting planes \cite{helmberg2000spectral}, and
bundle methods with generalized stabilization \cite{kiwiel1999bundle, frangioni2002generalized,
benamor2009on, frangioni2014generalized, de2016doubly}. 
Bundle methods are also widely used in non-convex 
optimization~\cite{schramm1992version, kiwiel1996restricted, luksan1998abm,
fuduli2004minimizing, haarala2007globally}.

Incremental bundle 
methods~\cite{emiel2010incremental} are based on a principle of 
selectively skipping oracle calls for some agents, while replacing them by an 
approximation.
Essentially, the function is evaluated when lower and upper estimates on function values 
are not sufficient~\cite{Ackooij2018IncrementalBM}.
This setting was analyzed within a framework of inexact oracles~\cite{emiel2010incremental, 
Oliveira2014ConvexPB, Oliveira2015ABM, van2016inexact}, 
\ie, noisy oracles that are ``asymptotically exact''.
Further, there has been an interest in asynchronous bundle methods as
well~\cite{iutzeler2020async, fischer2022asynchronous},
where the agent oracle calls have varying finite running times.

Our method for parameter discovery relies on a standard level bundle 
method~\cite{lemarechal1995new} for a few 
initial iterations. We then set the value of the proximal 
parameter to a Lagrange multiplier
and proceed with the proximal bundle method.
In a separate study, \cite{de2016doubly} propose an algorithm that automatically 
chooses between proximal and level bundle approaches at every iteration.
Another related line of research is based on variable metric
bundle methods~\cite{kiwiel1990proximity, lemarechal1995new, 
Lemarchal1997VariableMB, kiwiel2000efficiency, rey2002dynamical, 
frangioni2002generalized, Oliveira2014ConvexPB, Ackooij2018IncrementalBM}. 
In general, these approaches are not crafted for solving composite function 
minimization satisfying our access conditions.

Our setting with difficult-to-evaluate components has been 
considered in bundle methods before~\cite{emiel2010incremental}.
In those studies, the authors use the inexact evaluations to skip the ``hard'' 
components. 
We, on the other hand, focus on exact oracle evaluations or queries for each 
component $f_i$ alongside the structured function $g$,  
which is known as partially exact bundle method~\cite{van2017probabilistic}.
Our structured function $g$ presents a way to incorporate constraints into
bundle methods, which has also been considered before~\cite{kiwiel1990proximity}.

In our method we disaggregate the minorant of $f$, an idea which
has been extensively studied in the prior literature~\cite{frangioni2002generalized, 
lemarechal2009bundle, frangioni2014generalized, frangioni2014bundle, 
frangioni2020standard}.
When combined with the structured function $g$, which is not approximated by a minorant
but handled exactly, we arrive at what is referred to as
a disaggregate partially exact bundle method.

\paragraph{Software packages.}
Despite the large literature on related methods, relatively few 
open-source software packages are available.
We found only one open-source implementation
of a bundle method that is compatible with our access requirements.
\texttt{BundleMethod.jl}~\cite{BundleMethod.jl.0.3.2} is a Julia package
with implementations of proximal bundle method~\cite{kiwiel1990proximity}
and trust region bundle method~\cite{kim2019asynchronous}.
We found that in all our examples OSBDO works substantially better. 

We also mention here other existing open-source implementations of bundle 
methods that do not meet our access methods.
A Fortran implementation~\cite{makela2003multiobjective} with Julia interface
for multi-objective proximal bundle method~\cite{makela2003multiobjective, makela2016proximal},
available at~\cite{Karmitsa2016},
requires the objective function to have full domain.
A Fortran implementation~\cite{teo2010bundle},
uses bundle methods for unconstrained regularized risk minimization.
A limited memory bundle method~\cite{karmitsa2007lmbm, haarala2007globally,karmitsa2010limited} 
is available as a 
Fortran implementation at \cite{karmitsa2007lmbm},
solves nonsmooth large-scale minimization problems, either
unconstrained~\cite{karmitsa2007lmbm, haarala2007globally} or
bound constrained~\cite{karmitsa2010limited}.
An unconstrained proximal bundle method~\cite{diaz2021optimal}
with Julia implementation is available at~\cite{diaz2021}.

\subsection{Outline}
We describe our assembled bundle method for oracle-structured distributed optimization
in~\S\ref{s-bundle-method}.
In~\S\ref{s-agents} we take a deeper look at the agent functions $f_i$,
which often involve additional variables that are optimized
by the agent, and not used by the algorithm, which we call private 
variables.
We present several numerical examples in~\S\ref{s-examples}.
A convergence proof for our algorithm is given in~\S\ref{s-convergence}.

\section{Disaggregate partially exact bundle method} \label{s-bundle-method}
Here we describe our assembled bundle method to solve the oracle-structured 
distributed optimization problem \eqref{e-prob}.
We use the superscript $k$ to denote a vector or function at iteration 
$k$, as in
$x_i^k$ or $x^k=(x_1^k,\ldots, x_m^k)$, our estimates of 
$x_i^\star$ and $x^\star$ at iteration $k$.
Each iteration involves querying the agent objective functions, \ie,
evaluating $f_i(\tilde x_i^k)$ and a subgradient $q_i^k \in \partial
f_i(\tilde x_i^k)$ at a query point $\tilde x_i^k$, for $i=1,\ldots ,k$,
along with some computation that updates the iterates $x^k$.
(We will describe what the specific query points are later.)

\subsection{Minorants}
The basic idea in a bundle or cutting-plane method is to maintain and refine
a minorant of each agent function, denoted $\hat f_i^k:\reals^{n_i}\to \reals$.
(Minorant means that these functions satisfy 
$\hat f_i^k(x_i) \leq f_i(x_i)$ for all $x_i \in \reals^{n_i}$.)
They are constructed from initial (given) minorants $\hat f_i^0$,
and evaluations of the value and subgradients in previous iterations, as
\BEQ\label{e-fi-minorant}
\hat{f}_i^{k+1}(x_i) = \max \left(\hat{f}_i^k(x_i),~
f_i(\tilde{x}_i^{k+1}) + (q_i^{k+1})^T (x_i - \tilde{x}_i^{k+1}) 
\right), \quad i=1,\ldots, M.
\EEQ
Here we use the basic subgradient inequality
\[
f_i(\tilde x_i^{k+1}) + (q_i^{k+1})^T (x_i - \tilde x_i^{k+1}) \leq f_i(x_i)
\]
for all $x_i \in \reals^{n_i}$, which shows that the lefthand side,
which is an affine function of $x_i$, is a minorant of $f_i$.
We also note that the minorant \eqref{e-fi-minorant} 
is tight at $\tilde x_i^{k+1}$,
\ie,
\[
\hat{f}_i^{k+1}(\tilde x_i^{k+1}) = f_i(\tilde x_i^{k+1}),
\quad i=1, \ldots, M.
\]
From these agent objective minorants we obtain a minorant of the oracle objective $f$,
\BEQ\label{e-f-minorant}
\hat{f}^{k+1}(x) = \hat{f}^{k+1}_1(x_1)+ \cdots+  \hat f^{k+1}_M(x_M),
\EEQ
and in turn, a minorant of $h$,
\BEQ\label{e-h-minorant}
\hat h^{k+1}(x) = \hat f^{k+1}(x)+g(x).
\EEQ
The minorant of $f$ in \eqref{e-f-minorant} is referred to as a 
disaggregated minorant, to distinguish it from forming a minorant directly of $f$.

\paragraph{Initial minorant.}
The simplest initial minorant $\hat f_i^0$ is a constant,
a known lower bound on the agent objective $f_i(x_i)$.  
Later we will see how more
sophisticated minorants for the agent objectives can be obtained.
With a simple constant initial minorant for each agent, 
the minorant $\hat f^k$ is also piecewise affine function, since it is a sum of
$m$ terms, each the maximum of $k$ affine functions.

\paragraph{Lower bound on optimal value.}
Minorants allow us to compute a lower bound on $h^\star$,
the optimal value of the problem \eqref{e-prob}.
Since $\hat h^k$ is a minorant of $h$, we have
\BEQ \label{e-Lk}
L^k = \min_x \hat h^k(x) \leq h^\star.
\EEQ
Evaluating $L^k$ involves solving an optimization problem,
minimizing $g$ plus the piecewise affine function $\hat f^k$. 
Of course $h(x^k)$ is an upper bound on $h^\star$, so we have
\[
L^k \leq h^\star \leq h(x^k).
\]

\paragraph{Gap-based stopping criterion.}
There is a multitude of stopping criteria that have been proposed for 
bundle methods~\cite{lemarechal1995new, Lemarchal1997VariableMB, lemarechal1996bundle, 
frangioni2020standard, lemarechal2001lagrangian, hiriart1996convex, hiriart2013convex}. 
Many of them are based on approximately satisfying the optimality conditions.
For example, we can terminate when the norm of the aggregate
subgradient and aggregate linearization error are both 
small~\cite[Ch.~XIV]{Lemarchal1997VariableMB, lemarechal1996bundle, frangioni2020standard, 
hiriart1996convex}. 
Stopping can also be done based on the predicted function 
decrease~\cite{hiriart1996convex, frangioni2020standard}, 
which differs depending on the type of stabilization
in the bundle method~\cite[Ch.~XV]{hiriart1996convex}.
We refer the reader to the 
books~\cite{hiriart1996convex, hiriart2013convex}
for more details on stopping criterion for bundle methods.
There are other reasonable choices of stopping criteria, for 
example based on the minorant error $h(x^k)-\hat h^k(\tilde x^{k+1})$
and a subgradient of $\hat h^k$ at $\tilde x^{k+1}$ 
see, \eg, \cite[Ch.~5.1]{lemarechal2001lagrangian}.

Since we are in the regime where the agents are expensive to evaluate,
solving a minimization problem~\eqref{e-Lk} is not an issue, and especially if it is
not done every iteration, \ie, 
we compute the lower bound $L^k$ every few iterations.
The lower bound \eqref{e-Lk} gives us a gap-based stopping criterion.
In particular, this stopping criterion was previously used
in level bundle methods~\cite[Ch.~XIV]{lemarechal1995new, hiriart1996convex}.
We stop if either the absolute gap is small,
\BEQ\label{e-stopping-crit}
h(x^k)-L^k \leq \epsilon^\text{abs},
\EEQ
where $\epsilon^\text{abs}>0$ is a given absolute gap tolerance,
or the relative gap is small,
\BEQ\label{e-stopping-crit-rel}
h(x^k)-L^k \leq \epsilon^\text{rel} 
\min \{ |h(x^k)|, |L^k| \} \quad \mbox{and} \quad h(x^k)L^k > 0,
\EEQ
where $\epsilon^\text{rel}>0$ is a given relative gap tolerance.
(The sign condition guarantees that $\min \{ |h(x^k)|, |L^k| \} \leq |h^\star|$.)
This guarantees that we terminate with a point $x^k$ with objective
value that is either within absolute difference 
$\epsilon^\text{abs}$, or 
within relative distance $\epsilon^\text{rel}$, of $h^\star$.

For later reference, we define the relative gap as 
\BEQ\label{e-rel-gap}
\omega^k = \left\{ \begin{array}{ll} \frac{h(x^k)-L^k}{\min\{|h(x^k)|, |L^k|\}} &
h(x^k) L^k >0\\
\infty & \mbox{otherwise}. \end{array} \right.
\EEQ
We define the true relative gap as 
\BEQ\label{e-true-rel-gap}
\omega_\text{true}^k = \frac{h(x^k)-h^\star}{|h^\star|},
\EEQ
for $h^\star \neq 0$.
The relative gap is an upper bound on the true relative gap,
\ie, we have $\omega^k \geq \omega_\text{true}^k$.  But $\omega^k$
is known at iteration $k$, whereas $\omega_\text{true}^k$ is not.

\subsection{Oracle-structured bundle method}
The basic bundle method for \eqref{e-prob} is given below.
It includes two parameters $\eta$ and $\rho$ (discussed later),
and is initialized with an initial guess $x^0$ of $x$, and the initial
minorants of the agent objective functions $\hat f_i^0$.

\begin{algdesc}{\sc Bundle method for oracle-structured distributed 
optimization}\label{alg-bundle}
    {\footnotesize
    \begin{tabbing}
    {\bf given} $x^0 \in \dom h$, $h(x^0)$,
initial minorants $\hat{f}_i^0$, 
    parameters $\eta \in 
    (0,1)$ and $\rho>0$.\\*[\smallskipamount]
    {\bf for $k=0,1,\ldots$} \\
    \qquad \= 1.\ \emph{Check stopping criterion.} Quit 
    if~\eqref{e-stopping-crit} or \eqref{e-stopping-crit-rel}
holds.\\
    \> 2.\ \emph{Tentative update.} 
$\tilde{x}^{k+1} = \argmin_{x}
\left (\hat{h}^k(x) + (\rho/2) \|x - x^k \|_2^2 \right )$.\\
\>\qquad\= Record 
$\hat h^k(\tilde x^{k+1})$ and $g(\tilde x^{k+1})$.\\
    \> 3. \emph{Query agents.}  Evaluate $f_i(\tilde{x}_i^{k+1})$ and 
$q_i^{k+1} \in \partial f_i(\tilde{x}_i^{k+1})$, for $i=1,\ldots,M$.\\
\> 4. Compute $h(\tilde
x^{k+1})=f_1(\tilde x_1^{k+1})+ \cdots + f_M(\tilde x_M^{k+1})+
g(\tilde x^{k+1})$.\\
    \> 5.  Compute
	 $\delta^k = h(x^k) - \left (\hat{h}^k\left(\tilde{x}^{k+1}\right)
      + (\rho/2) \left\| \tilde{x}^{k+1} - x^k \right\|_2^2 \right )$. \\
    \> 6. \emph{Update iterate.} Set 
    $x^{k+1} = \begin{cases} \tilde{x}^{k+1} 
    & h(x^k) -
        h\left(\tilde{x}^{k+1}\right)\geq \eta \delta^k \\
        x^k & \text{otherwise}.
        \end{cases}$ \\
    \> 7. \emph{Update minorants.}\\
\>\>  Update $\hat{f}_i^{k+1}$ using \eqref{e-fi-minorant},
for $i=1,\ldots, M$.\\
\>\> Update $\hat{f}^{k+1}$ and $\hat h^{k+1}$ using
 \eqref{e-f-minorant} and \eqref{e-h-minorant}. 
    \end{tabbing}}  
\end{algdesc}

\paragraph{Comments.}
In step~1 we evaluate $L^k$, which involves solving 
an optimization problem, \ie, minimizing $\hat h^k(x)$.
We recognize step~2 as evaluating the proximal operator of $\hat h^k$ at 
$x^k$ \cite{parikh2014proximal}.
This step also involves solving an optimization problem, minimizing $\hat h^k(x)$
plus the quadratic (proximal) term $(\rho/2)\|x-x^k\|_2^2$.
We note that we always have $\tilde x^{k+1} \in \dom g \subseteq \dom f$,
so $h(\tilde x^{k+1})$ is always finite.
Step~3, the agent query, can be done in parallel across all agents.
Steps~4 and~5 involve only simple arithmetic using quantities already computed in
steps~2 and~3.
Steps~1,~5, and~6 use the quantity $h(x^k)$,
but that has already been computed, since $x^k$ is equal to a previous
$\tilde x^j$ for some $j \leq k$, and so was computed in step~4 of 
a previous iteration.

The substantial computation in each iteration of the bundle method is
the evaluation of the lower bound in step~1, evaluating 
the proximal operator of $\hat h$ in step~2, and
querying each agent oracle in step~3.
The computation of $L^k$ in step~1 is only used for the stopping criterion, 
so this step can be carried out only every few steps, to reduce the 
average computational burden.

\paragraph{Descent method.}
The quantity $\delta^k$ computed in step~5 is nonnegative.
To see this, we note from step~2 that $x=\tilde x^{k+1}$ minimizes
$\hat h^k(x) + (\rho/2)\|x-x^k\|_2^2$, so
\BEQ
\label{e-delt-pos}
\begin{array}{rcl}
\hat h^k(\tilde x^{k+1}) + (\rho/2)\|\tilde x^{k+1}-x^k\|_2^2
&\leq&
\hat h^k(x^k) + (\rho/2)\|x^k-x^k\|_2^2\\
&=&
\hat h^k(x^k)\\
&=& h(x^k),
\end{array}
\EEQ
from which $\delta^k \geq 0$ follows.
From step~6 we see that the bundle method is a 
descent method, \ie, $h(x^{k+1}) \leq h(x^k)$.
More specifically, 
$h(x^{k+1}) <h(x^k)$ if the tentative step is accepted,
\ie, $x^{k+1}=\tilde x^{k+1}$, and
$h(x^{k+1})=h(x^k)$ if the tentative step is not accepted,
\ie, $x^{k+1}=x^k$.

\paragraph{Convergence.}
In appendix~\ref{s-convergence} for completeness we give a proof that
the bundle method converges, \ie, $h(x^k) \to h^\star$ as $k \to \infty$.
Convergence proofs for bundle methods have long history,
dating back to the 1970s~\cite{lemarechal1978nonsmooth}.

\paragraph{Choice of parameters.}
The bundle method is not particularly sensitive to the choice of $\eta$;
the value $\eta = 0.01$ works well.
The value of $\rho$, however, can have a strong influence on 
the practical performance of the algorithm.  
We discuss choices of $\rho$ later in \S\ref{s-rho-choice}.

\paragraph{Dual variable.}
An estimate of an optimal dual variable $q^\star \in \partial f(x^\star)$
can be found when we compute the lower bound $L^k$ in step~1.
To do this we compute $L^k$ by solving the modified problem
\[
\begin{array}{ll} \mbox{minimize} & \hat f^{k+1}(x) + g(\tilde x)\\
\mbox{subject to} & \tilde x=x,
\end{array}
\]
with variables $x\in \reals^n$ and $\tilde x \in \reals^n$.
(This is the consensus form; see, \eg, \cite[Chap.~7]{boyd2011distributed}.)
The optimal dual variable associated with the constraint is our estimate
of $q^\star$.

\subsection{Diagonal preconditioning}\label{s-precond}
Practical convergence of the bundle method is greatly enhanced by
diagonal preconditioning.
This means that we choose a diagonal matrix $D$ with positive entries,
and define the (scaled) variable $\overline x =D^{-1}x$.   
Then we solve the problem
\eqref{e-prob} with variable $\overline x$ and functions 
\[
\overline f(\overline x) = f(D \overline x), \qquad
\overline g(\overline x) = g(D \overline x).
\]
Note that $\overline g$ is structured, since $g$ is.
We recover the solution of the original problem \eqref{e-prob}
as $x^\star = D\overline x^\star$.

The idea of preconditioning has a long-standing history dating back to 1845~\cite{jacobi1845ueber} 
and is a crucial technique in optimization; see, \eg, 
\cite{hestenes1952methods,sinkhorn1964relationship,
concus1985block,nocedal1999numerical,bradley2010algorithms,takapoui2016preconditioning}.
It is also closely related to the idea
of variable metric methods, where essentially a different
(not necessarily diagonal) preconditioning is applied each step.
Variable metric bundle methods have been studied in
\cite{lemarechal1994approach, bacaud2001bundle, helmberg2017dynamic}.
In contrast to our computationally cheap preconditioning technique, 
which is computed once before the algorithm starts,
existing variable metric bundle methods are considerably more 
complex~\cite{kiwiel1990proximity, lemarechal1995new, 
Lemarchal1997VariableMB, kiwiel2000efficiency, rey2002dynamical, 
frangioni2002generalized, Oliveira2014ConvexPB, Ackooij2018IncrementalBM};
in addition, most are not compatible with our specific access conditions.

Given $\overline x$, we query agent $i$ and the point $x_i = D_i\overline x_i$,
where $D_i$ is the submatrix of $D$ associated with $x_i$.
Agent $i$ responds with $f_i(x_i)=\overline f_i(\overline x_i)$
and $q_i \in \partial f_i(x_i)$.  The algorithm uses the function value
without change, since $f_i(x_i) = \overline f_i(\overline x_i)$,
and the scaled subgradient
\[
\overline q_i = D q_i \in \partial \overline f_i(\overline x_i).
\]

Diagonal scaling can be thought of as a thin layer or interface 
between the algorithm, which works with the scaled variable $\overline x$ 
and functions $\overline f_i$
and $\overline g$, and the agents,
which work with the original variables $x_i$ and the original functions $f_i$.
In particular, neither the algorithm nor the agents need to know that diagonal
scaling is being used, provided the query points and returned subgradients
are scaled correctly.

\paragraph{A specific choice for diagonal scaling.}
The matrix $D$ scales the original variable.  Our goal
is to scale the variables so they range over similar intervals,
say of width on the order of one.
To do this, we assume that we have
known lower and upper bound on the entries of $x$,
\BEQ\label{e-vable-bnds}
l \leq x \leq u,
\EEQ
where $l < u$.
(Presumably these constraints are included in $g$.)
With these variable bounds we can choose
\BEQ \label{e-precond-D}
D = \diag(u-l).
\EEQ
The new variable bounds have the form $\overline l \leq z \leq \overline u$,
with $\overline u - \overline l = \ones$, the vector with all entries one.

\subsection{Proximal parameter discovery}\label{s-rho-choice}
While the bundle algorithm converges for any positive 
value of the parameter $\rho$, fast convergence in practice requires
a reasonable choice that is somewhat problem dependent;
see, \eg,~\cite{diaz2021optimal}.
Several schemes can be used to find a good value of $\rho$.
In one common scheme $\rho$ is updated in each iteration depending on
various quantities computed in the iteration.  (As an example of such an adaptive
scheme, see \cite[\S 3.4.1]{boyd2011distributed}.)
Another general scheme is to start with some steps that are meant to discover 
a good value of $\rho$, as in \cite{rey2002dynamic}. 
For both such schemes, we fix $\rho$ after some modest number $K$ of iterations, 
so our proof of convergence (which assumes a fixed value of $\rho$) still
applies.
Our experiments suggest that a natural $\rho$-discovery method, described
below,  works well in practice. 
The idea of switching between proximal and level bundle methods has been proposed 
in \cite{de2016doubly}.  This method is more complex than our proposed simple approach,
which switches just once, after a fixed number of iterations.

For the discovery steps, we modify the update in step~2, where we minimize
\BEQ\label{e-prox-step}
\hat h^k(x) + (\rho/2) \|x-x^k\|_2^2,
\EEQ
\ie, evaluate the proximal operator of $\hat h^k$ at $x^k$, to solving the closely
related problem
\BEQ \label{e-prob-proj-subl}
\begin{array}{ll} \mbox{minimize} & (1/2)\|x - x^k \|_2^2\\
\mbox{subject to} & \hat h^k(x) \leq \eta^k,
\end{array}
\EEQ
where $\eta^k \in (L^k,h(x^k))$.
(This ensures that the problem is feasible, and that the constraint is tight.)
The problem \eqref{e-prob-proj-subl} finds the projection of the point $x^k$ 
onto the $\eta^k$-sublevel set of the minorant.
The idea of projecting onto the sublevel set instead of 
carrying out an explicit proximal step can be found in, \eg,~\cite[Ch.~2.3]{frangioni2020standard}.

It is easy to show that any solution of \eqref{e-prob-proj-subl} is also a solution
of \eqref{e-prox-step}, for some value of $\rho$.  Indeed we can find this value as 
$\rho = 1/\lambda$, where $\lambda$ is an optimal dual variable for the 
constraint in \eqref{e-prob-proj-subl}.
In other words: When we solve \eqref{e-prob-proj-subl}, we are actually computing
the proximal operator, \ie, solving \eqref{e-prox-step}, for a value of 
$\rho$ that we only find after solving it.

Our $\rho$-discovery method uses the update \eqref{e-prob-proj-subl} for the first
$K=20$ steps, with 
\[
\eta^k = \frac{h(x^k)+L^k}{2}.
\]
After $K$ steps, we use the standard update \eqref{e-prox-step} 
with the value of $\rho$ chosen as the geometric mean of the 
last $5$ values found during the discovery steps.

\subsection{Finite memory}
\label{finite-memory}
Evidently the optimization problems that must be solved to 
compute $L^k$ in step~1 and $\tilde x^{k+1}$ in step~2
grow in size as $k$ increases.
A standard method in bundle or cutting-plane type methods is to use 
a finite memory or constraint dropping version, also known as
bundle compression~\cite{Correa1993ConvergenceOS,  Oliveira2015ABM}.
The essential element that underpins theoretical 
convergence in compressed bundle methods
is aggregate linearization~\cite{kiwiel1983aggregate, Correa1993ConvergenceOS}.
It is given by the linearization of the minorant,
\[
l_i^{k+1}(x_i) = \hat{f}_i^{k}(\tilde x_i^{k+1}) + 
(\hat q_i^{k+1})^T(x_i - \tilde x_i^{k+1}),
\]
where $\hat q_i^{k+1} \in \partial \hat{f}_i^{k}(\tilde x_i^{k+1})$.
As a result, a finite memory version replaces 
the minorant \eqref{e-fi-minorant} with
\[
\hat{f}_i^{k+1}(x_i) = \max \left( l_i^{k+1}(x_i), 
\max_{j= \max\{0,k-m+2\}, \ldots, k} 
\left( f_i(\tilde{x}_i^{j+1}) + (q_i^{j+1})^T (x_i - \tilde{x}_i^{j+1}) 
\right)
\right),
\]
for all $i=1,\ldots, M$.
Thus we use only the last $m-1$ subgradients and values in the minorant,
instead of all previous ones with an additional affine term for aggregate 
linearization. This results in total of $m$ affine functions.
For further information on bundle compression, refer to \cite{Correa1993ConvergenceOS, 
lemarechal2001lagrangian, Oliveira2015ABM, frangioni2020standard}.
Interestingly, the minorant $\hat{f}_i^{k+1}(x_i)$ can be reduced to just 
two affine pieces. However, in practical applications, a smaller bundle size 
leads to slower convergence rates~\cite{Oliveira2015ABM}.
(Finite-memory should not be confused with limited-memory,
which can refer to a quasi-Newton algorithm that approximates
the Hessian using a finite number of function and gradient evaluations.)

With finite memory it is possible that the lower bounds
$L^k$ are not monotone nondecreasing.  In this case we can
keep track of the best (\ie, largest) lower bound found so far,
for use in our stopping criterion.

\section{Agents}\label{s-agents}
In this section we discuss some details of the agent objective functions,
and how to compute the value and a subgradient.  For simplicity we 
drop the subscript $i$ that was used denote agent $i$, so in this section,
we denote $x_i$ as $x$, $f_i$ as $f$, $q_i$ as $q$, and so on.

\paragraph{Private variables and partial minimization.}
In many cases the agent function $f$
is defined via partial minimization, as the optimal value of 
a problem with variable $x$ and additional variables $z$.
Specifically, $f(x)$ is the optimal value of the problem
\[
\begin{array}{ll}
\mbox{minimize} & F_0(x,z) \\ 
\mbox{subject to} & F_i(x,z)\leq 0,\quad i=1,\ldots,m, \\
& H_i (x,z)=0,\quad i=1,\ldots,p,
\end{array}
\]
with variable $z$.
(In this optimization problem, $x$ is a parameter.)
Here $F_i$ are jointly convex in $(x,z)$, and $H_i$ are jointly affine in $(x,z)$.
This is called partial minimization \cite[\S 4.1]{boyd2004convex}, and defines 
$f$ that is convex.
To evaluate $f$ we must solve a convex optimization problem.
We refer to $x$ as the public variable for the agent, and $z$ as its 
private variable, since its value (or even its existence) is not known
outside the agent.

We now explain how to find a subgradient $q\in \partial f(x)$.
We solve the equivalent convex problem
\[
\begin{array}{ll}
\mbox{minimize} & F_0(\tilde x,z) \\ 
\mbox{subject to} & F_i(\tilde x,z)\leq 0,\quad i=1,\ldots,m, \\
& H_i (\tilde x,z)=0,\quad i=1,\ldots,p,\\
& \tilde x=x,
\end{array}
\]
with variables $z$ and $\tilde x$.  (As in the problem above,
$x$ is a parameter in this optimization problem.)
We assume that strong duality holds for this problem 
(which is guaranteed if the stronger form of Slater's condition holds),
with $\nu$ denoting an optimal dual variable for the constraint $\tilde x=x$.
Then it is easy to show that $q = -\nu$ is a subgradient of 
$f$ \cite{364b-notes} assuming strong duality holds. 

\paragraph{Soft constraints and slack variables.}
We assume that the domain of the agent objective function includes
the domain of $g$.
This is critical since the agents will always be queried at a point $\tilde x
\in \dom g$, and we need the agent objective function to be finite for
any such $\tilde x$.
In some cases this property  does not hold for the natural definition 
of the agent objective function.
Here we explain how to modify an original definition of an agent 
objective function $f$ so that it does.

Let $\tilde f$ be the original agent objective function, 
which does not satisfy $\dom \tilde f \subseteq \dom g$.
It is often possible to replace constraints that appear in the original 
definition of $\tilde f$ with soft constraints, that ensure that
$\dom \tilde f \subseteq \dom g$ holds.

There is also a generic method that uses slack variables 
to ensure that $f$ has full domain (which implies the domain condition).
We define
\[
f(x) = \min_{\tilde x} \left( \tilde f(\tilde x) + \lambda \|\tilde x-x\|_1  
\right),
\]
where $\lambda >0$ is a parameter.  This $f$ is convex and has full domain.
For $\lambda$ large enough and $x\in \dom \tilde f$, we have $\tilde x =x$.
We can think of $\tilde x-x$ as a slack variable, used to guarantee that 
$f(x)$ is defined for all $x$.  We interpret $\lambda$ as a penalty for
using the slack variable.

\section{Examples} \label{s-examples}
In this section we present a number of examples to illustrate our method.
There are better methods to solve each of these examples,
which exploit the custom structure of the particular problem.
Our purpose here is to show that a variety of 
large-scale practical problems achieve good practical performance,
\ie, modest accuracy in some tens of iterations,
with OSBDO with default parameters.

We use the default parameters for the first four examples, except that 
we continue iterations after the algorithm would have terminated 
with the default 
tolerances $\epsilon^\text{abs}=10^{-3}$ and $\epsilon^\text{rel}=10^{-2}$,
to show the continued progress. 
The final subsection gives experimental results for the 
bundle method with finite memory
(described in \S\ref{finite-memory}).


\subsection{Supply chain}\label{s-supply-chain}
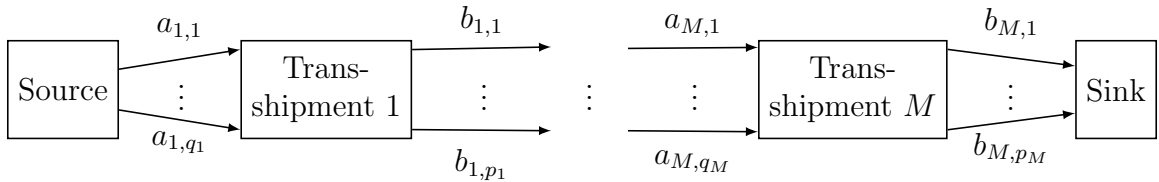
\begin{figure}
\centering
\begin {tikzpicture}[-latex ,auto ,node distance =1.0 cm and 3.5cm ,on grid ,
semithick ,
state/.style ={ rectangle ,top color =white , bottom color = white,
draw,black , text=black , minimum width = 1 cm, minimum height = 1.3cm}]
\node[state] (A) {Source};
\node[state, align=center] (C) [right=3.5cm of A]  {Trans-\\shipment $1$};
\node[state, draw=none] (E) [right=3.5cm of C]  {$\vdots$};
\node[state, align=center] (F) [right=3.5cm of E]  {Trans-\\shipment $M$};
\node[state] (H) [right=3.5cm of F]  {Sink};
\node[state, draw=none] [right=0.45*3.5cm of A]  {$\vdots$};
\node[state, draw=none]  [right=0.6*3.5cm of C]  {$\vdots$};
\node[state, draw=none] [right=0.4*3.5cm of E]  {$\vdots$};
\node[state, draw=none]  [right=0.6*3.5cm of F]  {$\vdots$};

\path (A.20) edge node[above=0.1cm] {$a_{1,1}$} (C.155);
\path (A.-20) edge node[below]{$a_{1,q_1}$} (C.205);
\path (C.25) edge node{$b_{1,1}$} (E.132);
\path (C.-25) edge node[below=0.1cm] {$b_{1,p_1}$} (E.227);
\path (E.47) edge node {$a_{M,1}$} (F.156);
\path (E.-47) edge node[below=0.1cm]{$a_{M,q_M}$} (F.204);
\path (F.23) edge node[above=0.1cm]{$b_{M,1}$} (H.150);
\path (F.-23) edge node[below]{$b_{M,p_M}$} (H.210);
\end{tikzpicture}
\caption{Supply chain consisting of a source (left), 
a sequence of $M$ trans-shipment components (middle), and a sink (right).}
\label{f-supply-chain}
\end{figure}

Supply chain problems involve the
placement and movement of inventory, including sourcing from a supplier
and distribution to an end customer.
A comprehensive review can be found in~\cite{trisna2016multi}.
Here we consider a single commodity network composed of 
a series of $M$ trans-shipment components.
Each of these has a vector
$a_i$ of (nonnegative) flows into it, and 
a vector $b_i$ of (nonnegative) flows out of it.
These are connected in series, with the first trans-shipment component's
output connected to the second component's input, and so on, so
\BEQ\label{e-flow-conserv}
b_1=a_2, ~\ldots, ~ b_{M-1} = a_M.
\EEQ
The vector $a_1$ gives the flows into the first trans-shipment component,
and $b_M$ is the vector of flows out of the last trans-shipment
component.
This is illustrated in figure~\ref{f-supply-chain}.

We assume the trans-shipment components are lossless, which means that
\BEQ\label{e-tranship-flow-balance}
\ones^T a_i=\ones^T b_i, \quad
i=1, \ldots, M,
\EEQ
\ie, the total flow of the commodity into
each trans-shipment component equals the total flow leaving it.

Each trans-shipment component has a nonnegative objective function $f_i(a_i,b_i)$ 
which we interpret as the cost of shipping the commodity from the input
flows to the output flows.
In addition there is a source (purchase) objective term $\psi^\text{src}(a_1)$,
the cost of purchasing the commodity, and 
a sink (delivery) objective term $\psi^\text{sink}(b_M)$,
which we interpret as the negative revenue derived from delivering the 
commodity. (So we expect that $\psi^\text{src}(a_1)$ is nonnegative, and
$\psi^\text{sink}(b_M)$ is nonpositive.)

The overall objective is the total of the trans-shipment costs and
the source and sink costs,
\[
\psi^\text{src}(a_1) + f_1(a_1,b_1) + \cdots + f_{M}(a_{M},b_{M}) +
\psi^\text{sink}(b_M).
\]
Our goal is to choose the input and output flows $a_i$ and $b_i$,
subject to the flow conservation constraints \eqref{e-flow-conserv},
so as to minimize this total objective.
We assume $\psi^\text{src}$, $\psi^\text{sink}$, and all $f_i$ are
convex, and that $\psi^\text{src}$ and $f_i$ are nonnegative.

\paragraph{Oracle-structured form.}
We put the supply chain problem into the form \eqref{e-prob} as follows.
We take $x_i=(a_i,b_i)$ for $i=1, \ldots, M$, with variable range
$0 \leq x \leq u$, where $u$ is an upper bound on the flows.
We take $g$ to be the source and sink objective terms, plus the
indicator function of the flow conservation constraints
\eqref{e-flow-conserv}, variable ranges, and flow balance
\eqref{e-tranship-flow-balance}:
\[
g(x) = \left\{ \begin{array}{ll} \psi^\text{src}(a_1)+\psi^\text{sink}(b_M) &
\eqref{e-flow-conserv}, ~ \eqref{e-tranship-flow-balance}, ~ 0 \leq x \leq u,\\
\infty & \mbox{otherwise}.
\end{array}\right.
\]
Roughly speaking, $g(x)$ is the gross negative profit, and $f(x)$ is
the total shipping cost.

Our initial minorant is the indicator function of the flow 
conservation constraint \eqref{e-tranship-flow-balance}, \ie,
$0$ if \eqref{e-tranship-flow-balance} holds and $\infty$ otherwise.
Note that this includes the lower bound $0\leq f_i(x_i)$.

\paragraph{Trans-shipment cost.}
The trans-shipment cost for agent $i$ is based on a fully bipartite
graph, with flow along each edge connecting one of $q_i$ inputs to 
each of $p_i$ outputs.
We represent the edge flows as $X_i \in \reals^{p_i \times q_i}$,
where $a_i \in \reals^{q_i}$ and $b_i\in \reals^{p_i}$, and $(X_i)_{jk}$
is the flow from input $k$ to output $j$.

We assume that each edge has a convex quadratic cost of the form
\[
(D_i)_{jk}(X_i)_{jk} + (E_i)_{jk} (X_i)_{jk}^2,
\]
with $(D_i)_{jk} \geq 0$, and in addition
is capacitated, \ie, $0\leq (X_i)_{jk}\leq (C_i)_{jk}$, with $(C_i)_{jk} \geq 0$.
Due to the capacity constraints, the domain condition
$\dom f_i \supseteq \dom g$ need not hold, so in addition we include a 
slack variable.
We define $f_i(x_i)$ as the optimal value of the trans-shipment problem
\[
\begin{array}{ll}
\mbox{minimize} & \sum_{j,k} \left((D_i)_{jk}(X_i)_{jk} + (E_i)_{jk} (X_i)
_{jk}^2\right) + 
\lambda \|r\|_1\\
\mbox{subject to} & 0 \leq (X_i)_{jk} \leq (C_i)_{jk}, \quad j=1,\ldots, p_i, \quad
k= 1, \ldots, q_i\\
& \sum_j (X_i)_{jk} =  (\tilde a_i)_k, \quad k=1, \ldots, q_i\\
& \sum_k (X_i)_{jk} =  (\tilde b_i)_j, \quad j=1, \ldots, p_i,\\
& (\tilde a_i,\tilde b_i)-r=x_i,
\end{array}
\]
with variables $(X_i)_{jk}$, $\tilde a_i$, $\tilde b_i$, and $r$.
Here $\lambda$ is a large positive parameter that penalizes using 
slack variables, \ie, having $\tilde a_i \neq a_i$ or $\tilde b_i \neq b_i$.

To evaluate $f_i(x_i)$ we solve the problem above, which is a 
convex QP.  To obtain a subgradient of $f_i$ at $x_i$
we use the negative optimal dual variables associated with the last constraint as $q_i$. 

\paragraph{Source and sink costs.}
We take simple linear source and sink costs,
\[
\psi^{src}(a_1) = \alpha^T a_1, \qquad
\psi^{sink}(b_M) = \beta^T b_M,
\]
with $\alpha \geq 0$ and $\beta \leq 0$.
We can interpret $\alpha_k$ as the price of acquiring the commodity at input
$k$ of the first trans-shipment component,
and $-\beta_k$ as the price of selling the commodity at output
$k$ of the last trans-shipment component.

\paragraph{Problem instance.}
We consider a problem instance with $M=5$ trans-shipment components, 
and dimensions of $a_i,b_i$ given by
\[
(20,30), \quad
(30,40), \quad
(40,25), \quad
(25,35), \quad
(35,20).
\]
We choose the edge capacities $(C_i)_{jk}$ from a log normal distribution,
with $\log (C_i)_{jk} \sim \mathcal N(0,1)$.
From these edge capacities we construct an upper bound on each component of 
$x_i$, as the maximum of the sum of the capacities of all edges that feed the
flow variable, and the sum of the capacities of all edges that flow out of it.
This gives us our upper bound $u$. 
The linear cost coefficients $(E_i)_{jk}$ are log normal,
with $\log (E_i)_{jk} \sim \mathcal N(0.07,0.7)$. 
The quadratic cost coefficients $(D_i)_{jk}$ are obtained from the capacity 
and linear cost coefficients as 
\[
(D_i)_{jk} = (E_i)_{jk}/ (2 (C_i)_{jk}).
\]
The sale prices $\alpha$ are uniform on $[8,10]$ and retail prices  
$-\beta$ are chosen uniformly on $[10, 12]$.
The total problem size is $n=300$ public variables (the input and output flows 
of the trans-shipment components), with an additional $4975$
private variables in the agents (the specific flows along all edges of the 
trans-shipment components).
The overall problem is a QP with $5275$ variables.

\begin{figure}
    \begin{center}
    \includegraphics[width=0.8\textwidth]{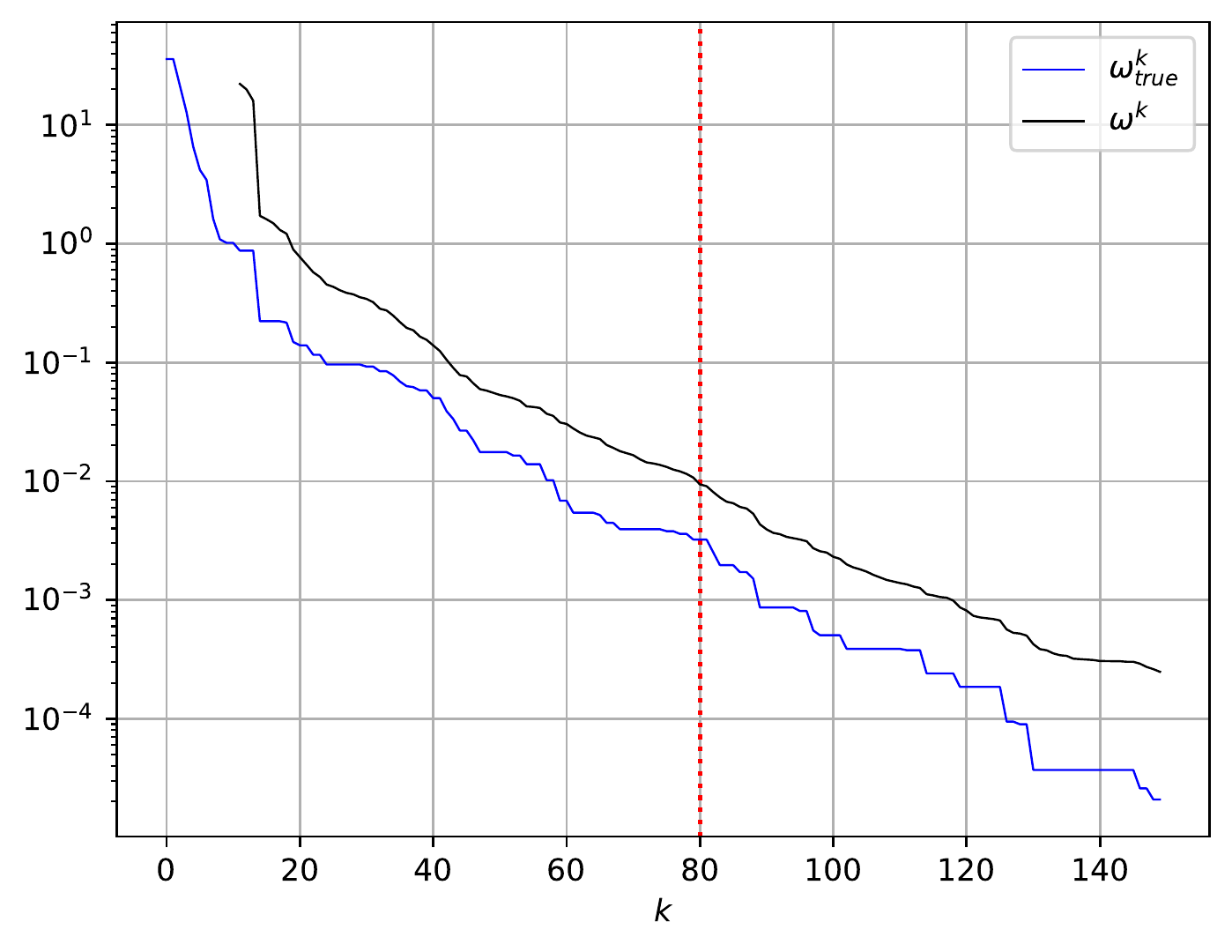}
    \end{center}
    \caption{Relative gap and true relative gap versus iterations 
for supply chain example.}
    \label{f-bundle}
\end{figure}

\paragraph{Results.}
Figure~\ref{f-bundle} shows the relative gap $\omega^k$
and true relative gap $\omega_\text{true}^k$, versus iterations.
With the default stopping criterion parameters
the algorithm would have terminated after $80$ iterations, when it can guarantee
that it is no more than $1\%$ suboptimal.   In fact, it is around $0.3\%$ suboptimal
at that point.  This is shown as the vertical dashed line in the plot.
We can also see that the relative accuracy was in fact better than $1\%$ after only
$59$ iterations, but, roughly speaking, we did not know it then.

\clearpage
\subsection{Resource allocation}\label{s-res-alloc}
Resource allocation problems consider how to allocate a 
limited amount of several resources to a number of participants 
in order to optimize some overall objective. 
This kind of problem arises in
communication networks \cite{han2008resource}, urban development 
\cite{wei2020simulation}, cloud computing \cite{choi2016optimization}, 
and many others. 
Various candidate algorithms have been explored to solve this 
problem, with interesting ones including ant colony algorithm 
\cite{yin2006ant}, genetic algorithm \cite{liu2005optimization}, and 
a graph-based approach \cite{zhou2021novel}. In this section, we 
demonstrate how bundle method can be exploited to address a distributed 
version of the generic resource allocation problem.

\paragraph{Resource allocation problem.}
We consider the optimal allocation of $n$ resources to $N$ participants.
We let $r_i \in \reals_+^n$ denote the amounts of the resources allocated
to participant $i$, for $i=1, \ldots, N$.
The utility derived by participant $i$ is $U_i(r_i)$, where $U_i: \reals_+^n \to
\reals$ is a concave nondecreasing utility function.
The resource allocation problem is to allocate
resources to maximize the total utility subject to a limit on the 
total resources allocated:
\[
\begin{array}{ll}
\mbox{maximize} & \sum_{i=1}^N U_i(r_i)  \\
\mbox{subject to} & r_i \geq 0, \quad i=1, \ldots, N\\
& \sum_{i=1}^N r_i \leq R,
\end{array}
\]
with variables $r_1, \ldots, r_N$,
where $R \in \reals_+^n$ is the total resources to be allocated, \ie, the budget.
We denote the optimal value, \ie, the maximum total utility, as a function of
$R$ as $U^\star(R)$.  It is also concave and nondecreasing.
When we solve this problem, an optimal dual variable associated with the 
last (budget) constraint can be interpreted as the prices of the resources.

\paragraph{Distributed resource allocation problem.} 
We have $M$ groups of participants, each with its own set of 
participants, resource budget $R_i \in \reals_+^n$, and utility $U^\star_i(R_i)$.
The distributed resource allocation problem is
\[
\begin{array}{ll}
\mbox{maximize} & \sum_{i=1}^{M} U^\star_i(R_i)  \\
\mbox{subject to} & R_i \geq 0, \quad i=1, \ldots, M\\
& \sum_{i=1}^M R_i \leq R,
\end{array}
\]
with variables $R_1, \ldots, R_M$,
where $R$ is the total budget of resources.
This problem has exactly the same form as the resource allocation problem, 
but here $U_i^\star(R_i)$
is the optimal total utility for group $i$ of participants, whereas in
the resource allocation problem, $U_i(r_i)$ is the utility of the single 
participant $i$.

\paragraph{Oracle-structured form.}
Each agent is associated with a group in the distributed resource allocation problem.
We take $x_i = R_i$, the total resource allocated to the participants 
in group $i$.
We take agent objective functions
\[
f_i(x_i) = -U_i^\star(x_i), \quad i=1, \ldots, M,
\]
the optimal (negative) total utility for its group of participants, given
resources $x_i$.
We take the structured objective function to be
\[
g(x) = \left\{ \begin{array}{ll} 0 &
x_1 + \cdots + x_M \leq R, \quad x_i \geq 0, \quad i=1, \ldots, M\\
\infty & \mbox{otherwise}.
\end{array} \right.
\]
With these agent and structured objectives, the problem \eqref{e-prob}
is equivalent to the distributed resource allocation problem.
The resource allocations to the individual participants within each group
are private variables; the public variables are the total resources
allocated to each group.

To evaluate $f_i(\tilde x_i)$, we solve the resource allocation problem for group $i$.
To find a subgradient $q\in \partial f_i(\tilde x_i)$, we take the 
negative of an optimal dual variable in the resource allocation problem, \ie, 
the negative of the optimal prices. To obtain a range on each variable, we use $0 \leq x_i \leq R$.

\paragraph{Problem instance.}
Our example uses participant utility functions of the form
\[
U_i(r_i)=  \geomean \left(A_ir_i+b_i\right),
\]
where $\geomean(u) = \left( \prod_{i=1}^p u_i\right)^{1/p}$ for 
$u \in \reals_+^p$ is the geometric mean function.
The entries of $A_i$ are nonnegative, so $U_i$ is concave and nondecreasing. 
We choose $A_i$ to be column sparse, with around
$\frac{n}{10}$ columns chosen at random to be nonzero.
The nonzero entries in these columns are chosen as uniform on $[0,1]$. 
We choose entries of $b_i$ to be uniform on $[0,\frac{n}{10}]$.  
The resource budgets $R_i$ are chosen from 
$\log R_i \sim \mathcal N(\log \frac{n}{10},1)$. 
The initial minorant is given by
\[
\hat{f}^0=\sum_{i=1}^{M}\sum_{j=1}^{N_i} -\geomean(A_{ij}R+b_{ij}),
\]
the negative of the utility if all agents were given the full budget
of resources.

For the specific instance we consider, we take $n=50$ resources 
and $M=50$ agents, each of which allocates resources to $N_i=10$ participants.
(As a single resource allocation problem we would have $50$ resources and 
$500$ participants.)  The utility functions use $p=5$, \ie, each is the 
geometric mean of $5$ affine functions.

\begin{figure}
    \begin{center}
    \includegraphics[width=0.8\textwidth]{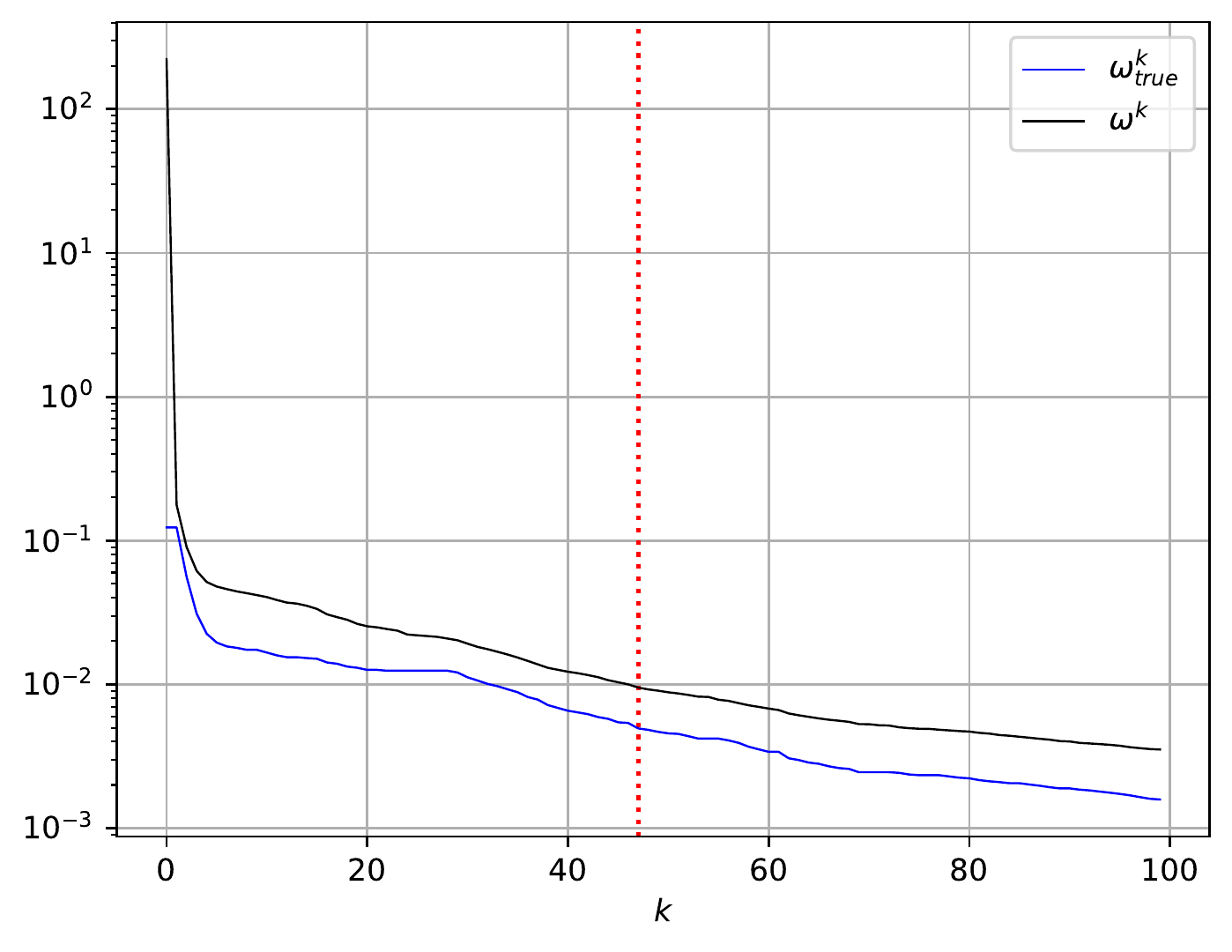}
    \end{center}
    \caption{Relative gap and true relative gap versus iterations 
for resource allocation example.}
    \label{ra-bundle}
\end{figure}

\paragraph{Results.}
Figure~\ref{ra-bundle} shows the relative gap $\omega^k$
and true relative gap $\omega_\text{true}^k$
versus iterations.
With the default stopping criterion the algorithm would have terminated
after $47$ iterations (shown as the vertical dashed line), when it can guarantee
that it is no more than $1\%$ suboptimal.   In fact, it is around $0.5\%$ suboptimal
at that point.
We can also see that the relative accuracy was better than $1\%$ after only
$33$ iterations.

\clearpage
\subsection{Multi-commodity flow}\label{s-mcf}
Multi-commodity flow problems involve shipping different commodities 
on the same network in a way such that the total utility is maximized,
while the total flow on each edge stays below its capacity. 
In~\cite{ouorou2000survey} 
authors give a survey of algorithms for this problem. 
We consider a network defined by a graph with $p$ vertices or nodes
and $q$ directed edges, defined by the incidence matrix 
$A\in \reals^{p \times q}$.
The network supports the flow of $M$ different commodities.
Each commodity has a source node, denoted $r_i \in \{1,\ldots, p\}$,
and a sink or destination node, denoted $s_i \in \{1,\ldots, p\}$.
The flow of commodity $i$ from the source to destination is given by $d_i\geq 0$.
We let $z_i\in \reals_+^q$ be the vector of flows of commodity $i$ 
on the edges.
Flow conservation is the constraint 
\[
Az_i +  d_i(e_{r_i} - e_{s_i})= 0, \quad i=1, \ldots, M,
\]
where $e_k$ denotes the $k$th unit vector, with $(e_k)_k=1$ and $(e_k)_j=0$ for
$j \neq k$.
Flow conservation requires that the flow is conserved at all nodes, 
with $d_i$ injected
at the source node, and $d_i$ removed at the sink node.
The utility of flow $i$ is $U_i(d_i)$,
where $U_i$ is a concave nondecreasing function.
Our objective is to maximize the total utility $U_1(d_1) + \cdots + U_M(d_M)$.

The total flow on the edges
must not exceed the capacities on the edges, given by $c\in \reals_+^q$, \ie,
\[
z_1+ \cdots+ z_M \leq c.
\]
(This capacity constraint couples the variables associated with 
the different commodities.) 

The variables in this multi-commodity flow problem are $z_1, \ldots, z_M$ and
$d_1, \ldots, d_M$.  
The data are the incidence matrix $A$, the edge capacities $c$,
the commodity source and sink nodes $(r_i, s_i)$, and the 
flow utility functions $U_i$.

It will be convenient to work with a form of the problem
where we split the capacity on each edge into $M$ different capacities for the 
different commodities.  We take
\[
z_i \leq c_i, \quad i=1, \ldots, M,
\]
where 
$c_1+ \cdots + c_M = c$ and $c_i \geq 0$ for $i=1, \ldots, M$.
We interpret $c_i$ as the edge capacity assigned to, or reserved for,
commodity $i$.
Our multi-commodity flow problem then has the form
\[
\begin{array}{ll} \mbox{maximize} & U_1(d_1) + \cdots + U_M(d_M)\\
\mbox{subject to} & 0 \leq z_i \leq c_i, \quad i=1, \ldots, M\\
& Az_i+d_i(e_{r_i} - e_{s_i}) =0, \quad i=1, \ldots, M\\
& c_1 + \cdots + c_M =c, \quad c_i \geq 0, \quad i=1, \ldots, M,
\end{array}
\]
with variables $z_i$, $d_i$, and $c_i$.  This is evidently equivalent 
to the original multi-commodity flow problem.

\paragraph{Oracle-structured form.}
We can put the multi-commodity flow problem into oracle-structured form
as follows.
We take $x_i=c_i$, the edge capacity assigned to commodity $i$.
We take the range of $x_i$ as $0\leq x_i \leq c$.
We take the agent cost function $f_i(x_i)$ to be the optimal value of
the single commodity flow problem (expressed as a minimization problem)
\[
\begin{array}{ll} \mbox{minimize} & - U_i(d_i)\\
\mbox{subject to} & 0 \leq z_i \leq x_i\\
& Az_i+d_i(e_{r_i} - e_{s_i}) =0,
\end{array}
\]
with variables $z_i$ and $d_i$. (Here $x_i=c_i$ is a parameter.)
These functions $f_i(x_i)$ are convex and nonincreasing in $x_i$, the 
capacity assigned to commodity $i$.
To evaluate $f_i$ we solve the single commodity flow problem above;
a subgradient $q_i$ is obtained as the negative optimal dual variable
associated with the capacity constraint $z_i \leq x_i$.

We take $g(x)$ to be the indicator function of 
\[
x_1 + \cdots + x_M =c, \quad x_i \geq 0, \quad i=1, \ldots, M
\]
(which includes the ranges of $x_i$).

All together there are $n=Mq$ variables, representing the allocation 
of edge capacity to the commodities.
There are also $M(q+1)$ private variables, which are the flows for each commodity
on each edge and the values of the flows of each commodity.

\paragraph{Problem instance.}
We consider an example with $M=10$ commodities,
and a graph with $p=100$ nodes and $q=1000$ edges. Edges are generated
randomly from pairs of nodes, with an additional cycle passing through 
all vertices (to ensure that the graph is strongly connected,
\ie, there is a directed path from any node to any other). We  
choose the source-destination pairs $(r_i, s_i)$ randomly.
We choose capacities $c_i$ from a uniform distribution on $[0.2,2]$.
The flow utilities $U_i$ are linear, \ie, $U_i(d_i)= b_i d_i$,
with $b_i$ chosen uniformly on $[0.5,1.5]$.
This problem instance has $n=10000$ variables, with an additional
$10010$ private variables.

\begin{figure}
\begin{center}
\includegraphics[width=0.8\textwidth]{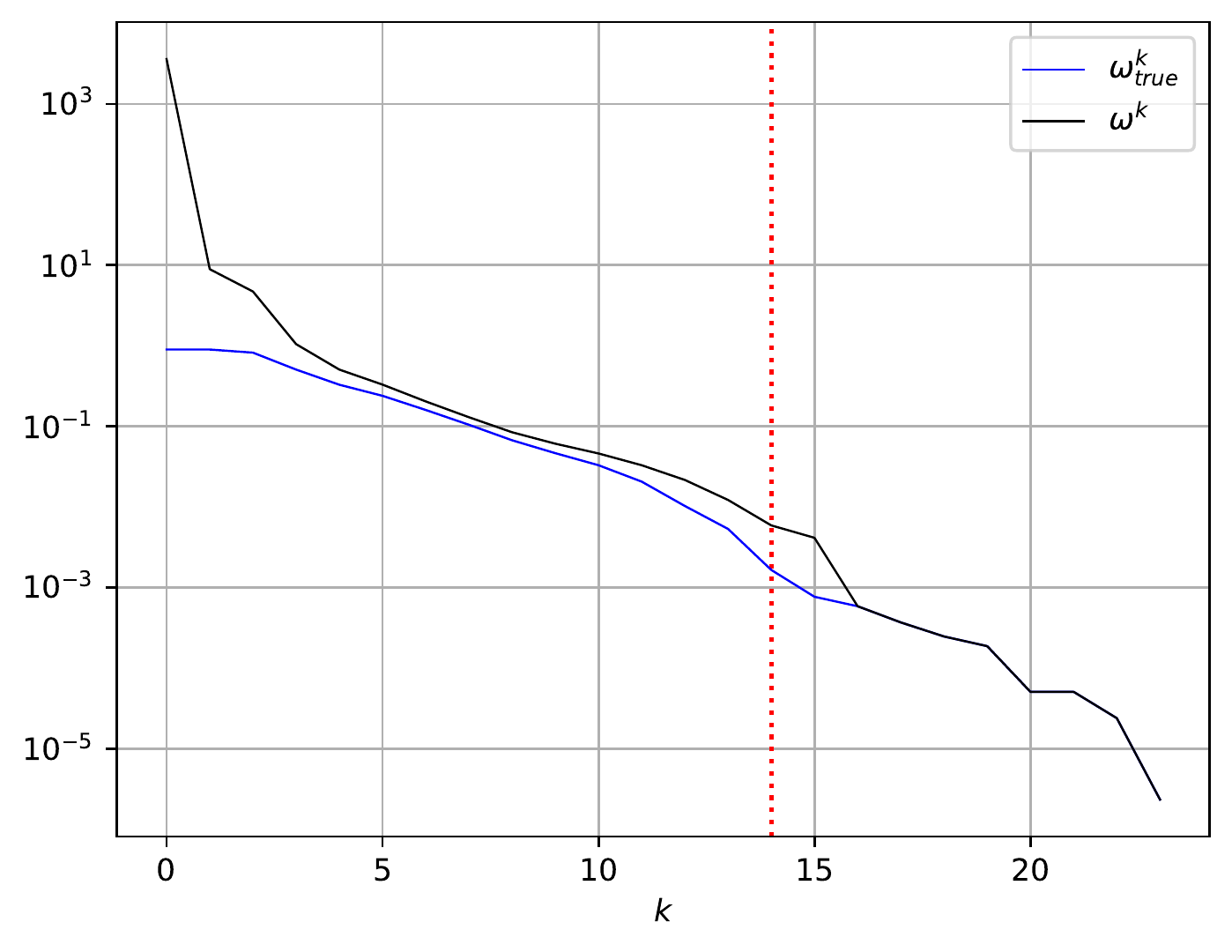}
\end{center}
\caption{Relative gap and true relative gap for multi-commodity flow example.}
\label{MCF_bundle}
\end{figure}

\paragraph{Results.}
Figure~\ref{MCF_bundle} shows the relative gap $\omega^k$
and true relative gap $\omega_\text{true}^k$
versus iterations.
With the default stopping criterion the algorithm would have terminated
after $14$ iterations (shown as the vertical dashed line), when it can guarantee
that it is no more than $1\%$ suboptimal.   In fact, it is around $0.6\%$ suboptimal
at that point.
We can also see that the relative accuracy was better than $1\%$ after only
$13$ iterations.

\clearpage
\subsection{Federated learning}\label{s-fed-learning}
Federated learning refers to distributed machine learning, where
agents keep their local data and collaboratively train a model using a distributed
algorithm. 
An overview of the development of federated learning is given in
\cite{li2020federated}.
In this section we consider the federated learning problem.

We are to fit a model parameter $\theta\in \reals^d$ 
to data that is stored in $M$ locations.
Associated with each location is a function $L_i:\reals^d \to \reals$,
where $L_i(\theta)$ is the loss for parameter value $\theta$ 
for the data held at location $i$.
We seek $\theta$ that minimizes
\[
\sum_{i=1}^M L_i(\theta)  + R(\theta),
\]
where $R: \reals^d \to \reals \cup \{\infty\}$ is a regularization 
function.
We assume that $L_i$ and $R$ are convex, so this fitting problem
is convex.
In federated learning~\cite{kairouz2021advances}, we solve the fitting problem
in a distributed manner, with each location handling its own data.

\paragraph{Oracle-structured form.}
We can put the federated learning problem into oracle-structured form
by taking $x_i$ to be the parameter estimate at location $i$,
$f_i = L_i$, and $g$ the indicator function for consensus plus the 
regularization,
\[
g(x) = \left\{ \begin{array}{ll}  R(x_1) & x_1 = \cdots = x_M\\
\infty & \mbox{otherwise}.
\end{array}\right.
\]

\paragraph{Problem instance.}
We consider a classification problem with logistic loss function,
\[
L_i(\theta) = \sum_{j=1}^{n_i} \log\left(1+\exp(-v_{ij} u_{ij}^T \theta)\right),
\]
where $v_{ij} \in \{-1,1\}$ is the label and $u_{ij} \in \reals^d$ is the feature
value for data point $j$ in location $i$, and $n_i$ is the number of data points
at location $i$.
We use $\ell_1$ regularization, \ie, $R(\theta)=\lambda \|\theta\|_1$, where 
$\lambda>0$.

Our example takes parameter dimension $d=500$, $M=10$ locations, 
and $n_i= 1000$ data points at each location.
In this problem there are no private variables, and the total dimension of 
$x$ is $n=Md=5000$.

We generate the data points as follows. The entries of $u_{ij}$ are 
$\mathcal N(0,1)$, and we take 
\[
v_{ij} = \sign \left(u_{ij}^T \theta^\mathrm{true} + z_{ij}\right),
\]
where $z_{ij} \sim \mathcal N(0,10^{-2})$ and $\theta^\mathrm{true}$ is 
a true value of the parameter, chosen as sparse with around $50$ nonzero 
entries, each $\mathcal N(0,1)$.
We choose $\lambda=5$.

\begin{figure}
\begin{center}
\includegraphics[width=0.8\textwidth]{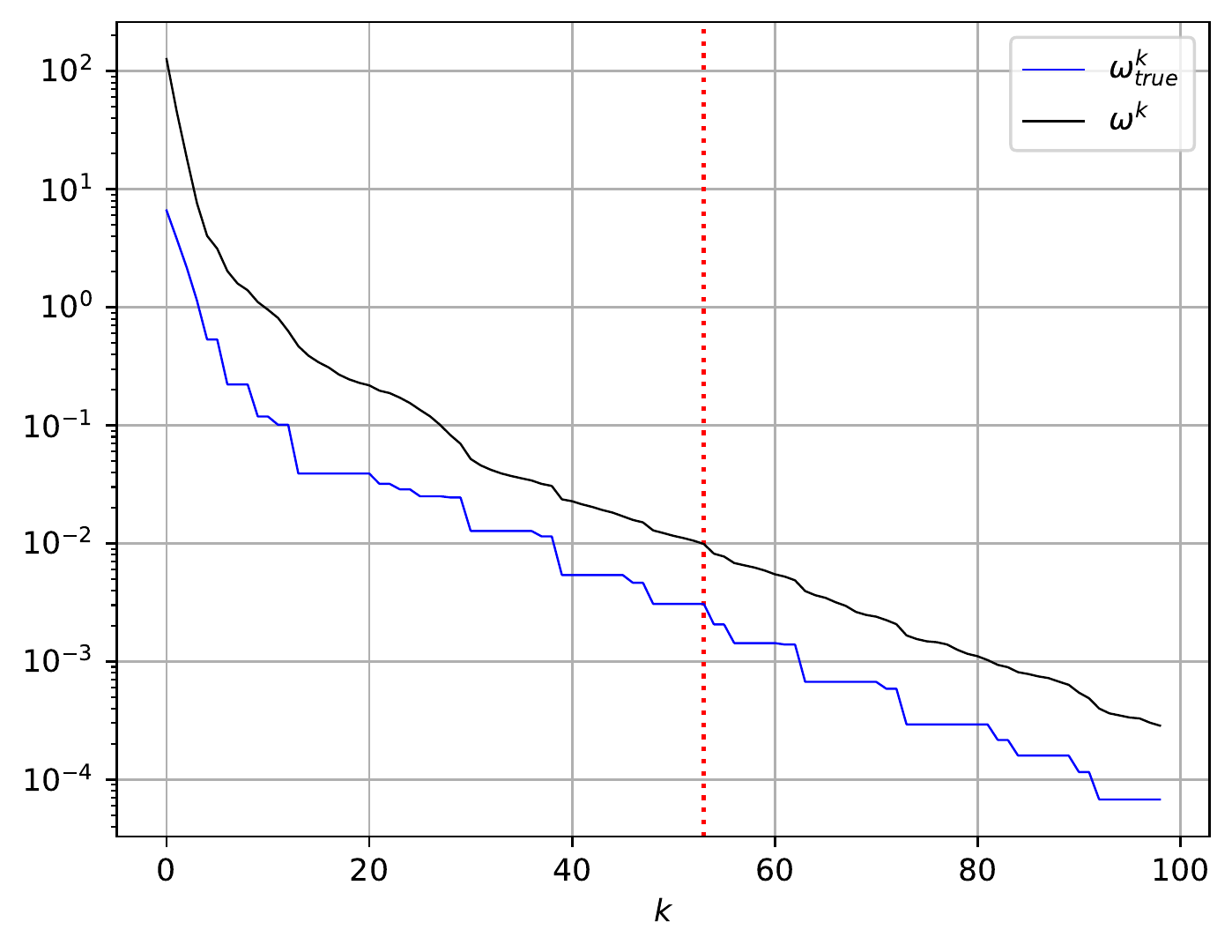}
\end{center}
\caption{Relative gap and true relative gap versus iterations
for federated learning example.}
\label{fl-bundle}
\end{figure}

\paragraph{Results.}
Figure~\ref{fl-bundle} shows the relative gap $\omega^k$
and true relative gap $\omega_\text{true}^k$
versus iterations.
With the default stopping criterion the algorithm would have terminated
after $53$ iterations (shown as the vertical dashed line), when it can guarantee
that it is no more than $1\%$ suboptimal.   In fact, it is around $0.3\%$ suboptimal
at that point.
We can also see that the relative accuracy was better than $1\%$ after only
$39$ iterations.

\clearpage
\subsection{Finite-memory experiments}
\begin{figure}
    \begin{center}
    \includegraphics[width=0.8\textwidth]{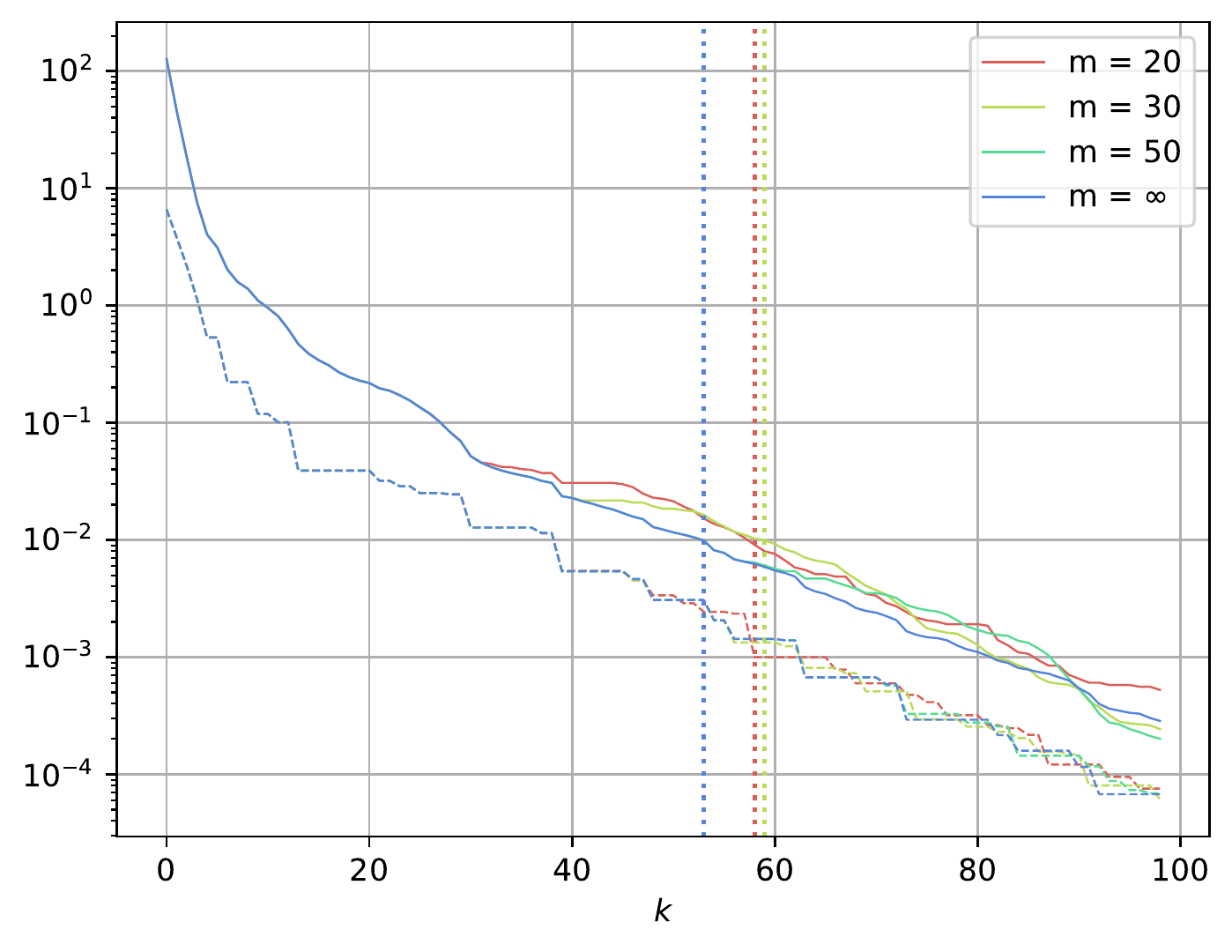}
    \end{center}
    \caption{Relative gap (solid) and true relative gap (dashed) 
versus iterations for federated learning example, 
with finite memory values $m=20$, $m=30$, $m=50$, and $m=\infty$.}
\label{f-FM-FL}
\end{figure}

In this subsection we present results showing how limiting the memory 
to various values affects convergence.  In many cases, limiting the value to
$m=20$ or more has negligible effect.  As an example, figure~\ref{f-FM-FL}
shows the effect on convergence of memory with values $m=20$, $m=30$, $m=50$,
and $m=\infty$ for the federated learning problem described above.
In this example finite memory has essentially small effect on the convergence.

\begin{figure}
    \begin{center}
    \includegraphics[width=0.8\textwidth]{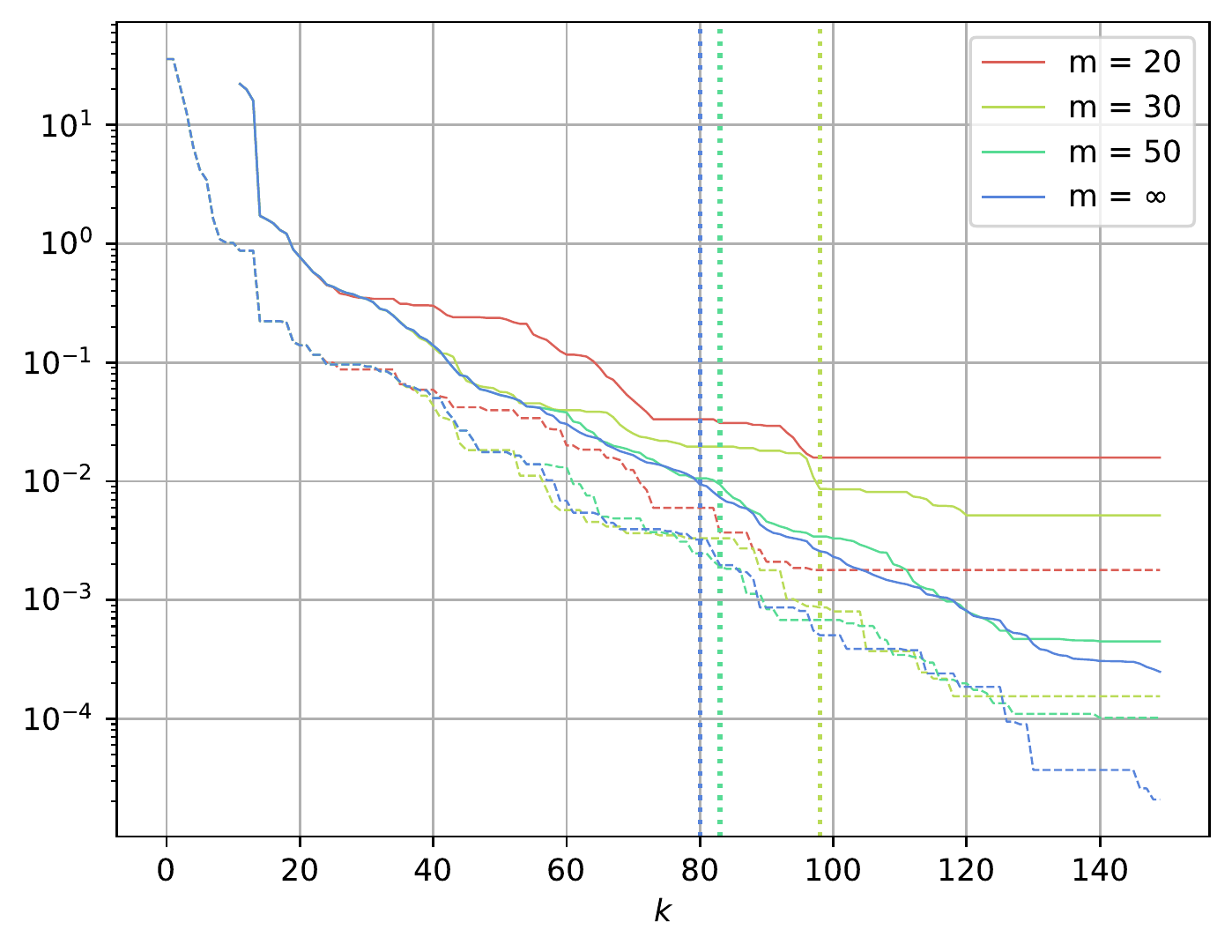}
    \end{center}
    \caption{Relative gap (solid) and true relative gap (dashed) 
versus iterations
for supply chain example, with finite memory values $m=20$, $m=30$, $m=50$, and $m=\infty$.}
\label{f-FM-SC}
\end{figure}

As an example of a case where finite memory does affect the convergence,
figure~\ref{f-FM-SC} shows the effect of finite memory with values $m=20$, $m=30$,
$m=50$, and $m=\infty$.
With $m=20$, the algorithm shows minimal improvement beyond double 
the number of iterations required for a method with $m=\infty$ to 
achieve the default tolerance level;
with $m=30$ there is a modest increase;
with $m=50$ there is a small increase.

\clearpage
\section{Conclusions}
We focus on developing a good practical method for distributed convex optimization
in a setting where the agents support a value/subgradient oracle, which can take
substantial effort to evaluate, and the coupling among the agent variables
is given explicitly and exactly as a structured convex problem, 
possibly including constraints.
(This differs from the more typical setting, where the agent functions
are differentiable and can be evaluated quickly, and the coupling function
has an analytical proximal operator.)
Our assumptions allow us to consider methods that carry out more computation
in each iteration, such as cutting-plane or bundle methods, that typically involve
the solution of a QP.
We have found that a basic bundle-type method, when combined with diagonal scaling
and a good algorithm parameter discovery method, gives good practical convergence across 
a number of problems types and sizes.  Here by good practical convergence we mean
that with default algorithm parameters, 
a reasonable approximate solution can be found in a few tens of iterations, 
and a higher accuracy solution (which is generally not needed in applications)
can be obtained in perhaps a hundred or fewer iterations.
(Theoretical convergence of the algorithm is always guaranteed.)

Our methods combines multiple variations of known techniques for bundle-type methods
into a solver has a number of attractive features.  
First, it 
has essentially no algorithm parameters, and works well with the few parameters
set to default values.  
Second, it achieves good 
practical convergence across a number of problems types and sizes.
Third, it can warm start when the coupling changes, by saving the 
information obtained in previous agent evaluations.

\section*{Acknowledgments}

We thank Parth Nobel, Nikhil Devanathan, Garrett van Ryzin, 
Dominique Perrault-Joncas, Lee Dicker, and Manan Chopra for very helpful
discussions about the problem and formulation.  
The supply chain example
was suggested by van Ryzin, Perrault-Joncas, and Dicker. The communication 
layer for the implementation with structured variables, to be described 
in a future paper, was designed by Parth Nobel and Manan Chopra.
We thank Mateo D{\'\i}az for pointing us to some very relevant literature
that we had missed in an early version of this paper. 
We thank three anonymous reviewers who gave extensive and helpful feedback on 
an early version of this paper.

We gratefully acknowledge support from Amazon, Stanford Graduate Fellowship, 
Office of Naval Research, and the Oliger Memorial Fellowship.
This research was partially supported by ACCESS –- AI Chip Center for Emerging
Smart Systems, sponsored by InnoHK funding, Hong Kong SAR.

\clearpage
\bibliography{references}

\newcommand{\etalchar}[1]{$^{#1}$}
\begin{thebibliography}{vABdOS17}

\bibitem[AF18]{Ackooij2018IncrementalBM}
W.~van Ackooij and A.~Frangioni.
\newblock Incremental bundle methods using upper models.
\newblock {\em SIAM Journal on Optimization}, 28:379--410, 2018.

\bibitem[AV95]{atkinson1995acp}
D.~S. Atkinson and P.~M. Vaidya.
\newblock A cutting plane algorithm for convex programming that uses analytic
  centers.
\newblock {\em Mathematical Programming}, 69:1--43, 1995.

\bibitem[AVDB18]{agrawal2018rewriting}
A.~Agrawal, R.~Verschueren, S.~Diamond, and S.~Boyd.
\newblock A rewriting system for convex optimization problems.
\newblock {\em Journal of Control and Decision}, 5(1):42--60, 2018.

\bibitem[BADF09]{benamor2009on}
H.~M.~T. Ben~Amor, J.~Desrosiers, and A.~Frangioni.
\newblock On the choice of explicit stabilizing terms in column generation.
\newblock {\em Discrete Applied Mathematics}, 157(6):1167--1184, 2009.

\bibitem[BDPV22]{364b-notes}
S.~Boyd, J.~Duchi, M.~Pilanci, and L.~Vandenberghe.
\newblock Stanford {EE} 364b, lecture notes: Subgradients, 2022.
\newblock URL:
  \url{https://web.stanford.edu/class/ee364b/lectures/subgradients_notes.pdf}.
  Last visited on 2022/08/04.

\bibitem[Bel05]{belloni2005lecture}
A.~Belloni.
\newblock Lecture notes for {IAP} 2005 course introduction to bundle methods,
  2005.

\bibitem[BLRS01]{bacaud2001bundle}
L.~Bacaud, C.~Lemar{\'e}chal, A.~Renaud, and C.~Sagastiz{\'a}bal.
\newblock Bundle methods in stochastic optimal power management: A
  disaggregated approach using preconditioners.
\newblock {\em Computational Optimization and Applications}, 20:227--244, 2001.

\bibitem[BMLRT15]{burachik2015anas}
R.~S. Burachik, J.~E. Mart{\'i}nez-Legaz, M.~Rezaie, and M.~Th{\'e}ra.
\newblock An additive subfamily of enlargements of a maximally monotone
  operator.
\newblock {\em Set-Valued and Variational Analysis}, 23:643--665, 2015.

\bibitem[BMR03]{birgin2003inexact}
E.~G. Birgin, J.~M. Mart{\'\i}nez, and M.~Raydan.
\newblock Inexact spectral projected gradient methods on convex sets.
\newblock {\em IMA Journal of Numerical Analysis}, 23(4):539--559, 10 2003.

\bibitem[BPC11]{boyd2011distributed}
S.~Boyd, N.~Parikh, and E.~Chu.
\newblock {\em Distributed Optimization and Statistical Learning via the
  Alternating Direction Method of Multipliers}.
\newblock Now Publishers, 2011.

\bibitem[BQ00]{burke2000ots}
J.~V. Burke and M.~Qian.
\newblock On the superlinear convergence of the variable metric proximal point
  algorithm using {B}royden and {BFGS} matrix secant updating.
\newblock {\em Mathematical Programming}, 88:157--181, 2000.

\bibitem[Bra10]{bradley2010algorithms}
A.~Bradley.
\newblock {\em Algorithms for the equilibration of matrices and their
  application to limited-memory Quasi-Newton methods}.
\newblock PhD thesis, Stanford University, CA, 2010.

\bibitem[Bru75]{bruck1975iterative}
R.~E. Bruck, Jr.
\newblock An iterative solution of a variational inequality for certain
  monotone operators in {H}ilbert space.
\newblock {\em Bulletin of the American Mathematical Society}, 81:890--892,
  1975.

\bibitem[BV04]{boyd2004convex}
S.~Boyd and L.~Vandenberghe.
\newblock {\em Convex Optimization}.
\newblock Cambridge University Press, 2004.

\bibitem[CF99]{chen1999proximal}
X.~Chen and M.~Fukushima.
\newblock Proximal quasi-{N}ewton methods for nondifferentiable convex
  optimization.
\newblock {\em Mathematical Programming}, 85(2):313--334, 1999.

\bibitem[CG59]{cheney1959newton}
E.~W. Cheney and A.~A. Goldstein.
\newblock Newton's method for convex programming and {T}chebycheff
  approximation.
\newblock {\em Numerische Mathematik}, 1:253--268, 1959.

\bibitem[CGM85]{concus1985block}
P.~Concus, G.~Golub, and G.~Meurant.
\newblock Block preconditioning for the conjugate gradient method.
\newblock {\em SIAM Journal on Scientific and Statistical Computing},
  6(1):220--252, 1985.

\bibitem[CL93]{Correa1993ConvergenceOS}
R.~Correa and C.~Lemar{\'e}chal.
\newblock Convergence of some algorithms for convex minimization.
\newblock {\em Mathematical Programming}, 62:261--275, 1993.

\bibitem[CL16]{choi2016optimization}
Y.~Choi and Y.~Lim.
\newblock Optimization approach for resource allocation on cloud computing for
  {IoT}.
\newblock {\em International Journal of Distributed Sensor Networks},
  12(3):3479247, 2016.

\bibitem[CP11]{combettes2011proximal}
P.~L. Combettes and J.-C. Pesquet.
\newblock Proximal splitting methods in signal processing.
\newblock In {\em Fixed-Point Algorithms for Inverse Problems in Science and
  Engineering}, pages 185--212. Springer, 2011.

\bibitem[CR97]{chen1997convergence}
G.~H. Chen and R.~T. Rockafellar.
\newblock Convergence rates in forward--backward splitting.
\newblock {\em SIAM Journal on Optimization}, 7(2):421--444, 1997.

\bibitem[DB16]{diamond2016cvxpy}
S.~Diamond and S.~Boyd.
\newblock {CVXPY}: {A} {P}ython-embedded modeling language for convex
  optimization.
\newblock {\em Journal of Machine Learning Research}, 17(83):1--5, 2016.

\bibitem[DG21]{diaz2021optimal}
M.~D{\'\i}az and B.~Grimmer.
\newblock Optimal convergence rates for the proximal bundle method, 2021.

\bibitem[DHS11]{duchi2011adaptive}
J.~Duchi, E.~Hazan, and Y.~Singer.
\newblock Adaptive subgradient methods for online learning and stochastic
  optimization.
\newblock {\em Journal of Machine Learning Research}, 12(7):2121--2159, 2011.

\bibitem[DV85]{dem1985nondifferentiable}
V.~F. Dem'yanov and L.~V. Vasil'ev.
\newblock {\em Nondifferentiable Optimization}.
\newblock Translations Series in Mathematics and Engineering. Springer New
  York, 1985.

\bibitem[Dí21]{diaz2021}
M.~Díaz.
\newblock proximal-bundle-method.
\newblock \url{https://github.com/mateodd25/proximal-bundle-method}, 2021.

\bibitem[EM75]{elzinga1975central}
J.~Elzinga and T.~G. Moore.
\newblock A central cutting plane algorithm for the convex programming problem.
\newblock {\em Mathematical Programming}, 8:134--145, 1975.

\bibitem[ES10]{emiel2010incremental}
G.~Emiel and C.~Sagastiz{\'a}bal.
\newblock Incremental-like bundle methods with application to energy planning.
\newblock {\em Computational Optimization and Applications}, 46(2):305--332,
  2010.

\bibitem[FG14a]{frangioni2014bundle}
A.~Frangioni and E.~Gorgone.
\newblock Bundle methods for sum-functions with ``easy'' components:
  Applications to multicommodity network design.
\newblock {\em Mathematical Programming}, 145:133--161, 2014.

\bibitem[FG14b]{frangioni2014generalized}
A.~Frangioni and E.~Gorgone.
\newblock Generalized bundle methods for sum-functions with ``easy''
  components: Applications to multicommodity network design.
\newblock {\em Mathematical Programming}, 145:133--161, 2014.

\bibitem[FGG04]{fuduli2004minimizing}
A.~Fuduli, M.~Gaudioso, and G.~Giallombardo.
\newblock Minimizing nonconvex nonsmooth functions via cutting planes and
  proximity control.
\newblock {\em SIAM Journal on Optimization}, 14(3):743--756, 2004.

\bibitem[Fis22]{fischer2022asynchronous}
F.~Fischer.
\newblock An asynchronous proximal bundle method.
\newblock {\em Optimization Online}, 2022.

\bibitem[Fra02]{frangioni2002generalized}
A.~Frangioni.
\newblock Generalized bundle methods.
\newblock {\em SIAM Journal on Optimization}, 13(1):117--156, 2002.

\bibitem[Fra20]{frangioni2020standard}
A.~Frangioni.
\newblock Standard bundle methods: {U}ntrusted models and duality.
\newblock In {\em Numerical Nonsmooth Optimization}, pages 61--116. Springer,
  2020.

\bibitem[GBY06]{grant2006disciplined}
M.~Grant, S.~Boyd, and Y.~Ye.
\newblock Disciplined convex programming.
\newblock In {\em Global Optimization}, pages 155--210. Springer, 2006.

\bibitem[GP79]{gonzaga1979constraint}
C.~Gonzaga and E.~Polak.
\newblock On constraint dropping schemes and optimality functions for a class
  of outer approximations algorithms.
\newblock {\em SIAM Journal on Control and Optimization}, 17(4):477--493, 1979.

\bibitem[Hin01]{hintermuller2001proximal}
M.~Hinterm{\"u}ller.
\newblock A proximal bundle method based on approximate subgradients.
\newblock {\em Computational Optimization and Applications}, 20(3):245--266,
  2001.

\bibitem[HL08]{han2008resource}
Z.~Han and K.~J.~R. Liu.
\newblock {\em Resource Allocation for Wireless Networks: Basics, Techniques,
  and Applications}.
\newblock Cambridge University Press, 2008.

\bibitem[HMM04]{haarala2004new}
M.~Haarala, K.~Miettinen, and M.~M. M{\"a}kel{\"a}.
\newblock New limited memory bundle method for large-scale nonsmooth
  optimization.
\newblock {\em Optimization Methods and Software}, 19(6):673--692, 2004.

\bibitem[HMM07]{haarala2007globally}
N.~Haarala, K.~Miettinen, and M.~M. M{\"a}kel{\"a}.
\newblock Globally convergent limited memory bundle method for large-scale
  nonsmooth optimization.
\newblock {\em Mathematical Programming}, 109:181--205, 2007.

\bibitem[HP17]{helmberg2017dynamic}
C.~Helmberg and A.~Pichler.
\newblock {\em Dynamic Scaling and Submodel Selection in Bundle Methods for
  Convex Optimization}.
\newblock Technische Universit{\"a}t Chemnitz, Fakult{\"a}t f{\"u}r Mathematik,
  2017.

\bibitem[HR00]{helmberg2000spectral}
C.~Helmberg and F.~Rendl.
\newblock A spectral bundle method for semidefinite programming.
\newblock {\em SIAM Journal on Optimization}, 10(3):673--696, 2000.

\bibitem[HS{\etalchar{+}}52]{hestenes1952methods}
M.~Hestenes, E.~Stiefel, et~al.
\newblock Methods of conjugate gradients for solving linear systems.
\newblock {\em Journal of Research of the National Bureau of Standards},
  49(6):409--436, 1952.

\bibitem[HSS16]{hare2016proximal}
W.~Hare, C.~Sagastiz{\'a}bal, and M.~Solodov.
\newblock A proximal bundle method for nonsmooth nonconvex functions with
  inexact information.
\newblock {\em Computational Optimization and Applications}, 63(1):1--28, 2016.

\bibitem[HUL96]{hiriart1996convex}
J.-B. Hiriart-Urruty and C.~Lemar{\'e}chal.
\newblock {\em Convex Analysis and Minimization Algorithms II: Advanced Theory
  and Bundle Methods}.
\newblock Grundlehren der mathematischen Wissenschaften. Springer Berlin
  Heidelberg, 1996.

\bibitem[HUL13]{hiriart2013convex}
J.-B. Hiriart-Urruty and C.~Lemar{\'e}chal.
\newblock {\em Convex Analysis and Minimization Algorithms I: Fundamentals},
  volume 305.
\newblock Springer Science \& Business Media, 2013.

\bibitem[IMdO20]{iutzeler2020async}
F.~Iutzeler, J.~Malick, and W.~de~Oliveira.
\newblock Asynchronous level bundle methods.
\newblock {\em {Mathematical Programming}}, 184:319--348, 2020.

\bibitem[Jac45]{jacobi1845ueber}
C.~Jacobi.
\newblock Ueber eine neue aufl{\"o}sungsart der bei der methode der kleinsten
  quadrate vorkommenden line{\"a}ren gleichungen.
\newblock {\em Astronomische Nachrichten}, 22(20):297--306, 1845.

\bibitem[Kar07]{karmitsa2007lmbm}
N.~Karmitsa.
\newblock {LMBM}--{FORTRAN} subroutines for large-scale nonsmooth minimization:
  User's manual'.
\newblock {\em TUCS Techn. Rep.}, 77:856, 2007.

\bibitem[Kar16]{Karmitsa2016}
N.~Karmitsa.
\newblock Proximal bundle method.
\newblock \url{http://napsu.karmitsa.fi/proxbundle/}, 2016.

\bibitem[Kel60]{kelley1960cutting}
J.~E. Kelley, Jr.
\newblock The cutting-plane method for solving convex programs.
\newblock {\em Journal of the Society for Industrial and Applied Mathematics},
  8(4):703--712, 1960.

\bibitem[Kiw83]{kiwiel1983aggregate}
K.~Kiwiel.
\newblock An aggregate subgradient method for nonsmooth convex minimization.
\newblock {\em Mathematical Programming}, 27:320--341, 1983.

\bibitem[Kiw85]{kiwiel1985algorithm}
K.~C. Kiwiel.
\newblock An algorithm for nonsmooth convex minimization with errors.
\newblock {\em Mathematics of Computation}, 45(171):173--180, 1985.

\bibitem[Kiw90]{kiwiel1990proximity}
K.~Kiwiel.
\newblock Proximity control in bundle methods for convex nondifferentiable
  minimization.
\newblock {\em Mathematical Programming}, 46(1-3):105--122, 1990.

\bibitem[Kiw95]{kiwiel1995approximations}
K.~C. Kiwiel.
\newblock Approximations in proximal bundle methods and decomposition of convex
  programs.
\newblock {\em Journal of Optimization Theory and applications},
  84(3):529--548, 1995.

\bibitem[Kiw96]{kiwiel1996restricted}
K.~C. Kiwiel.
\newblock Restricted step and {L}evenberg–{M}arquardt techniques in proximal
  bundle methods for nonconvex nondifferentiable optimization.
\newblock {\em SIAM Journal on Optimization}, 6(1):227--249, 1996.

\bibitem[Kiw99]{kiwiel1999bundle}
K.~C. Kiwiel.
\newblock A bundle {B}regman proximal method for convex nondifferentiable
  minimization.
\newblock {\em Mathematical Programming}, 85(2):241--258, 1999.

\bibitem[Kiw00]{kiwiel2000efficiency}
K.~C. Kiwiel.
\newblock Efficiency of proximal bundle methods.
\newblock {\em Journal of Optimization Theory and Applications},
  104(3):589--603, 2000.

\bibitem[Kiw06]{kiwiel2006proximal}
K.~C. Kiwiel.
\newblock A proximal bundle method with approximate subgradient linearizations.
\newblock {\em SIAM Journal on Optimization}, 16(4):1007--1023, 2006.

\bibitem[KM10]{karmitsa2010limited}
N.~Karmitsa and M.~M{\"a}kel{\"a}.
\newblock Limited memory bundle method for large bound constrained nonsmooth
  optimization: Convergence analysis.
\newblock {\em Optimization Methods \& Software}, 25(6):895--916, 2010.

\bibitem[KMA{\etalchar{+}}21]{kairouz2021advances}
P.~Kairouz, H.~B. McMahan, B.~Avent, A.~Bellet, M.~Bennis, A.~N. Bhagoji,
  K.~Bonawitz, Z.~Charles, G.~Cormode, R.~Cummings, et~al.
\newblock Advances and open problems in federated learning.
\newblock {\em Foundations and Trends{\textregistered} in Machine Learning},
  14(1--2):1--210, 2021.

\bibitem[KPZ19]{kim2019asynchronous}
K.~Kim, C.~Petra, and V.~Zavala.
\newblock An asynchronous bundle-trust-region method for dual decomposition of
  stochastic mixed-integer programming.
\newblock {\em SIAM Journal on Optimization}, 29(1):318--342, 2019.

\bibitem[KZNS21]{BundleMethod.jl.0.3.2}
K.~Kim, W.~Zhang, H.~Nakao, and M.~Schanen.
\newblock {BundleMethod.jl: Implementation of Bundle Methods in Julia}, March
  2021.

\bibitem[Lem75]{lemarechal1975aneo}
C.~Lemar{\'e}chal.
\newblock An extension of {D}avidon methods to non differentiable problems.
\newblock {\em Mathematical Programming Study}, 3:95--109, 1975.

\bibitem[Lem78]{lemarechal1978nonsmooth}
C.~Lemarechal.
\newblock Nonsmooth optimization and descent methods.
\newblock 1978.

\bibitem[Lem01]{lemarechal2001lagrangian}
C.~Lemar{\'e}chal.
\newblock Lagrangian relaxation.
\newblock In {\em Computational Combinatorial Optimization}, pages 112--156.
  Springer, 2001.

\bibitem[LM79]{lions1979Splitting}
P.~L. Lions and B.~Mercier.
\newblock Splitting algorithms for the sum of two nonlinear operators.
\newblock {\em SIAM Journal on Numerical Analysis}, 16(6):964--979, 1979.

\bibitem[LNN95]{lemarechal1995new}
C.~Lemar{\'e}chal, A.~Nemirovskii, and Y.~Nesterov.
\newblock New variants of bundle methods.
\newblock {\em Mathematical Programming}, 69(1):111--147, 1995.

\bibitem[LOP09]{lemarechal2009bundle}
C.~Lemar{\'e}chal, A.~Ouorou, and G.~Petrou.
\newblock A bundle-type algorithm for routing in telecommunication data
  networks.
\newblock {\em Computational Optimization and Applications}, 44:385--409, 2009.

\bibitem[LPM18]{lv2018proximal}
J.~Lv, L.~Pang, and F.~Meng.
\newblock A proximal bundle method for constrained nonsmooth nonconvex
  optimization with inexact information.
\newblock {\em Journal of Global Optimization}, 70(3):517--549, 2018.

\bibitem[LS94]{lemarechal1994approach}
C.~Lemar{\'e}chal and C.~Sagastiz{\'a}bal.
\newblock An approach to variable metric bundle methods.
\newblock In {\em System Modelling and Optimization}, pages 144--162. Springer,
  1994.

\bibitem[LS97]{Lemarchal1997VariableMB}
C.~Lemar{\'e}chal and C.~Sagastiz{\'a}bal.
\newblock Variable metric bundle methods: From conceptual to implementable
  forms.
\newblock {\em Mathematical Programming}, 76:393--410, 1997.

\bibitem[LSPR96]{lemarechal1996bundle}
C.~Lemar{\'e}chal, C.~Sagastiz{\'a}bal, F.~Pellegrino, and A.~Renaud.
\newblock Bundle methods applied to the unit-commitment problem.
\newblock In {\em System Modelling and Optimization: Proceedings of the
  Seventeenth IFIP TC7 Conference on System Modelling and Optimization, 1995},
  pages 395--402. Springer, 1996.

\bibitem[LSTS20]{li2020federated}
T.~Li, A.~K. Sahu, A.~Talwalkar, and V.~Smith.
\newblock Federated learning: Challenges, methods, and future directions.
\newblock {\em IEEE Signal Processing Magazine}, 37(3):50--60, 2020.

\bibitem[LV98]{luksan1998abm}
L.~Luk{\v s}an and J.~Vl{\v c}ek.
\newblock A bundle-{N}ewton method for nonsmooth unconstrained minimization.
\newblock {\em Mathematical Programming}, 83:373--391, 1998.

\bibitem[LV99]{luksan1999globally}
L.~Luk{\v s}an and J.~Vl{\v c}ek.
\newblock Globally convergent variable metric method for convex nonsmooth
  unconstrained minimization.
\newblock {\em Journal of Optimization Theory and Applications}, 102:593--613,
  1999.

\bibitem[LZDL05]{liu2005optimization}
Y.~Liu, S.~Zhao, X.~Du, and S.~Li.
\newblock Optimization of resource allocation in construction using genetic
  algorithms.
\newblock In {\em 2005 International Conference on Machine Learning and
  Cybernetics}, volume~6, pages 3428--3432. IEEE, 2005.

\bibitem[M{\"a}k03]{makela2003multiobjective}
M.~M{\"a}kel{\"a}.
\newblock Multiobjective proximal bundle method for nonconvex nonsmooth
  optimization: Fortran subroutine {MPBNGC 2.0}.
\newblock {\em Reports of the Department of Mathematical Information
  Technology, Series B. Scientific Computing, B}, 13:2003, 2003.

\bibitem[MHB75]{10.2307/169692}
R.~Marsten, W.~Hogan, and J.~Blankenship.
\newblock The boxstep method for large-scale optimization.
\newblock {\em Operations Research}, 23(3):389--405, 1975.

\bibitem[Mif77]{mifflin1977semismooth}
R.~Mifflin.
\newblock Semismooth and semiconvex functions in constrained optimization.
\newblock {\em SIAM Journal on Control and Optimization}, 15(6):959--972, 1977.

\bibitem[Mif96]{mifflin1996quasi}
R.~Mifflin.
\newblock A quasi-second-order proximal bundle algorithm.
\newblock {\em Mathematical Programming}, 73(1):51--72, 1996.

\bibitem[MKW16]{makela2016proximal}
M.~M{\"a}kel{\"a}, N.~Karmitsa, and O.~Wilppu.
\newblock Proximal bundle method for nonsmooth and nonconvex multiobjective
  optimization.
\newblock {\em Mathematical Modeling and Optimization of Complex Structures},
  pages 191--204, 2016.

\bibitem[Nes83]{nesterov1983amf}
Y.~Nesterov.
\newblock A method for solving the convex programming problem with convergence
  rate $\mathcal{O}(1/k^2)$.
\newblock {\em Proceedings of the USSR Academy of Sciences}, 269:543--547,
  1983.

\bibitem[NW99]{nocedal1999numerical}
J.~Nocedal and S.~Wright.
\newblock {\em Numerical Optimization}.
\newblock Springer, 1999.

\bibitem[OE15]{Oliveira2015ABM}
W.~de Oliveira and J.~Eckstein.
\newblock A bundle method for exploiting additive structure in difficult
  optimization problems.
\newblock {\em Optimization Online}, 2015.

\bibitem[OMV00]{ouorou2000survey}
A.~Ouorou, P.~Mahey, and J.-Ph. Vial.
\newblock A survey of algorithms for convex multicommodity flow problems.
\newblock {\em Management Science}, 46(1):126--147, 2000.

\bibitem[OS16]{de2016doubly}
W.~de Oliveira and M.~Solodov.
\newblock A doubly stabilized bundle method for nonsmooth convex optimization.
\newblock {\em Mathematical Programming}, 156(1):125--159, 2016.

\bibitem[OS20]{oliveira2020bundle}
W.~de Oliveira and M.~Solodov.
\newblock Bundle methods for inexact data.
\newblock In {\em Numerical Nonsmooth Optimization}, pages 417--459. Springer,
  2020.

\bibitem[OSL14]{Oliveira2014ConvexPB}
W.~de Oliveira, C.~Sagastiz{\'a}bal, and C.~Lemar{\'e}chal.
\newblock Convex proximal bundle methods in depth: a unified analysis for
  inexact oracles.
\newblock {\em Mathematical Programming}, 148:241--277, 2014.

\bibitem[Pas79]{passty1979ergodic}
G.~B. Passty.
\newblock Ergodic convergence to a zero of the sum of monotone operators in
  {H}ilbert space.
\newblock {\em Journal of Mathematical Analysis and Applications},
  72(2):383--390, 1979.

\bibitem[PB{\etalchar{+}}14]{parikh2014proximal}
N.~Parikh, S.~Boyd, et~al.
\newblock Proximal algorithms.
\newblock {\em Foundations and Trends{\textregistered} in Optimization},
  1(3):127--239, 2014.

\bibitem[Roc81]{rockafellar1981TheTO}
R.~T. Rockafellar.
\newblock {\em The Theory of Subgradients and its Applications to Problems of
  Optimization}.
\newblock Heldermann Verlag, 1981.

\bibitem[RS02a]{rey2002dynamical}
P.~Rey and C.~Sagastiz{\'a}bal.
\newblock Dynamical adjustment of the prox-parameter in bundle methods.
\newblock {\em Optimization}, 51(2):423--447, 2002.

\bibitem[RS02b]{rey2002dynamic}
P.~A. Rey and C.~Sagastizábal.
\newblock Dynamical adjustment of the prox-parameter in bundle methods.
\newblock {\em Optimization}, 51(2):423--447, 2002.

\bibitem[Sch22]{schechtman2022stochastic}
S.~Schechtman.
\newblock Stochastic proximal subgradient descent oscillates in the vicinity of
  its accumulation set.
\newblock {\em Optimization Letters}, pages 1--14, 2022.

\bibitem[Sho12]{shor2012minimization}
N.~Z. Shor.
\newblock {\em Minimization Methods for Non-differentiable Functions},
  volume~3.
\newblock Springer Science \& Business Media, 2012.

\bibitem[Sin64]{sinkhorn1964relationship}
R.~Sinkhorn.
\newblock A relationship between arbitrary positive matrices and doubly
  stochastic matrices.
\newblock {\em The Annals of Mathematical Statistics}, 35(2):876--879, 1964.

\bibitem[SNW12]{sra2012optimization}
Suvrit Sra, S.~Nowozin, and S.~J. Wright.
\newblock {\em Optimization for Machine Learning}.
\newblock Mit Press, 2012.

\bibitem[SZ92]{schramm1992version}
H.~Schramm and J.~Zowe.
\newblock A version of the bundle idea for minimizing a nonsmooth function:
  Conceptual idea, convergence analysis, numerical results.
\newblock {\em SIAM Journal on Optimization}, 2(1):121--152, 1992.

\bibitem[TJ16]{takapoui2016preconditioning}
R.~Takapoui and H.~Javadi.
\newblock Preconditioning via diagonal scaling.
\newblock 2016.

\bibitem[TMAS16]{trisna2016multi}
T.~Trisna, M.~Marimin, Y.~Arkeman, and T.~Sunarti.
\newblock Multi-objective optimization for supply chain management problem: A
  literature review.
\newblock {\em Decision Science Letters}, 5(2):283--316, 2016.

\bibitem[TVSL10]{teo2010bundle}
C.~H. Teo, S.V.N. Vishwanathan, A.~Smola, and Q.~Le.
\newblock Bundle methods for regularized risk minimization.
\newblock {\em Journal of Machine Learning Research}, 11(1), 2010.

\bibitem[vABdOS17]{van2017probabilistic}
W.~van Ackooij, V.~Berge, W.~de~Oliveira, and C.~Sagastiz{\'a}bal.
\newblock Probabilistic optimization via approximate $p$-efficient points and
  bundle methods.
\newblock {\em Computers \& Operations Research}, 77:177--193, 2017.

\bibitem[vAFdO16]{van2016inexact}
W.~van Ackooij, A.~Frangioni, and W.~de~Oliveira.
\newblock Inexact stabilized {B}enders’ decomposition approaches with
  application to chance-constrained problems with finite support.
\newblock {\em Computational Optimization and Applications}, 65:637--669, 2016.

\bibitem[WP95]{westerlund1995extended}
T.~Westerlund and F.~Pettersson.
\newblock An extended cutting plane method for solving convex {MINLP} problems.
\newblock {\em Computers $\&$ Chemical Engineering}, 19:131--136, 1995.

\bibitem[WZX{\etalchar{+}}20]{wei2020simulation}
F.~Wei, X.~Zhang, J.~Xu, J.~Bing, and G.~Pan.
\newblock Simulation of water resource allocation for sustainable urban
  development: An integrated optimization approach.
\newblock {\em Journal of Cleaner Production}, 273:122537, 2020.

\bibitem[YW06]{yin2006ant}
P.~Yin and J.~Wang.
\newblock Ant colony optimization for the nonlinear resource allocation
  problem.
\newblock {\em Applied Mathematics and Computation}, 174(2):1438--1453, 2006.

\bibitem[ZBL{\etalchar{+}}21]{zhou2021novel}
B.~Zhou, J.~Bao, J.~Li, Y.~Lu, T.~Liu, and Q.~Zhang.
\newblock A novel knowledge graph-based optimization approach for resource
  allocation in discrete manufacturing workshops.
\newblock {\em Robotics and Computer-Integrated Manufacturing}, 71:102160,
  2021.

\end{thebibliography}
\clearpage

\appendix
\section{Convergence proof} \label{s-convergence}
In this section we give a proof of convergence of the bundle method for 
oracle-structured optimization.  Our proof uses well known ideas, and borrows
heavily from~\cite{belloni2005lecture}.  
We will make one additional (and traditional) assumption, that $f$ and $g$ 
are Lipschitz continuous on $\dom g$.

We say that the update was accepted in iteration $k$ 
if $x^{k+1}=\tilde x^{k+1}$.
Suppose this occurs in iterations $k_1< k_2< \cdots$.
We let $K=\{k_1, k_2, \ldots \}$ denote the set of iterations where the update
was accepted.
We distinguish two cases: $|K|  = \infty$ and $|K|  <  \infty$.

\paragraph{Infinite updates.}
We assume $|K|  =  \infty$.
First we establish that $\delta^{k_s}\to 0$ as $s\to \infty$.
Since $k=k_s$ is an accepted step, from step~6 of the algorithm we have
\[
\eta\delta^{k_s} \leq  
h(x^{k_s})  - h(x^{k_{s}+1}) =
h(x^{k_s})  - h(x^{k_{s+1}}).
\]
Summing this inequality from $s=1$ to $s=l$ and dividing by $\eta$ gives
\[
\sum_{s=1}^l\delta^{k_s} \leq 
\frac{h(x^{k_1}) - h(x^{k_{l+1}})}{\eta} \leq
\frac{h(x^0) - h^\star}{\eta},
\]
which implies that $\delta^{k_s}$ is summable, and so 
converges to zero as $s \to \infty$.

Since $\tilde x^{k_s+1}$ minimizes $\hat h^{k_s}(x)+(\rho/2)\|x-x^{k_s}\|_2^2$,
we have
\[
\partial \hat{h}^{k_s}\left(\tilde x^{k_s+1} \right) +
\rho( \tilde x^{k_s+1}-x^{k_s})  \ni 0.
\]
Using $\tilde x^{k_s+1} = x^{k_s+1} = x^{k_{s+1}}$, we have
\[
\rho(x^{k_s}-x^{k_{s+1}}) 
\in \partial \hat{h}^{k_s}\left(x^{k_{s+1}} \right).
\]
It follows that
\[
h^\star  = h(x^\star) \geq \hat h^{k_s}(x^\star) \geq
\hat h^{k_s}(x^{k_{s+1}}) + 
\rho(x^{k_s}-x^{k_{s+1}}) ^T (x^\star - x^{k_{s+1}}).
\]
We first rewrite this as
\BEAS
\frac{h^\star - \hat h^{k_s}(x^{k_{s+1}})}{\rho}  &\geq&  
(x^{k_s}-x^{k_{s+1}}) ^T (x^\star - x^{k_s}) +
(x^{k_s}-x^{k_{s+1}}) ^T (x^{k_s}- x^{k_{s+1}}))\\
&=& 
(x^{k_s}-x^{k_{s+1}}) ^T (x^\star - x^{k_s}) +
\|x^{k_s}-x^{k_{s+1}}\|_2^2,
\EEAS
and then in the form we will use below,
\[
2(x^{k_s}-x^{k_{s+1}}) ^T (x^\star - x^{k_s}) \leq (2/\rho)\left(
h^\star - \hat h^{k_s}(x^{k_{s+1}})\right) -  2\|x^{k_s}-x^{k_{s+1}}\|_2^2.
\]
Now we use a standard subgradient algorithm argument.  We have
\BEAS
\|x^{k_{s+1}} - x^\star \|_2^2 &=& \|x^{k_s} - x^\star
 \|_2^2 + \| x^{k_{s+1}} - x^{k_s}\|_2^2 + 2 (x^{k_s} -
 x^{k_{s+1}})^T (x^\star - x^{k_s}) \\
&\leq& \|x^{k_s} - x^\star \|_2^2+(2/\rho) \left ( h^\star -
\hat{h}^{k_s}(x^{k_{s+1}}) \right) - \| x^{k_{s+1}} - x^{k_s}\|_2^2
\\
&=& \|x^{k_s} - x^\star \|_2^2+(2/\rho) \left ( h^\star - h(x^{k_s}) 
+\delta^{k_s} \right).
\EEAS 
Summing this inequality from $s=1$ to $s=l$ and re-arranging yields
\BEAS
(2/\rho) \sum_{s=1}^l \left (h(x^{k_s}) - h^\star\right)  &\leq&
\|x^{k_1} - x^\star \|_2^2 - 
\|x^{k_{l+1}} - x^\star \|_2^2 + (2/\rho )\sum_{s=1}^l\delta^{k_s} \\
&\leq& \|x^{k_1} - x^\star \|_2^2  +2(h(x^0) - h^\star)/\eta \rho.
\EEAS
It follows that the nonnegative series $h(x^{k_s}) - h^\star$ is summable,
and therefore, $h(x^{k_s})\to h^\star$ as $s \to \infty$.

\paragraph{Finite updates.}
We assume $|K|  <  \infty$, with $p = \max K$ its largest entry. 
It follows that
for any $k>p$, we have $ h(x^k) -
h\left(\tilde{x}^{k+1}\right)<\eta \delta^k $. Note that
$x^k=x^p$ for all $k\geq p+1$. Moreover, using
\[
\|\tilde{x}^{k+2} - x^p \|_2^2 =  \|\tilde{x}^{k+2} -
\tilde{x}^{k+1} \|_2^2 + \left\|\tilde{x}^{k+1} - x^p\right \|_2^2 -
2\left( x^p - \tilde{x}^{k+1}\right)^T(\tilde{x}^{k+2}
-\tilde{x}^{k+1}) 
\] 
with $\rho(x^p - \tilde{x}^{k+1}) \in \partial
\hat{h}^k\left(\tilde{x}^{k+1} \right)$ and
$\hat{h}^{k+1}(\tilde{x}^{k+2})\geq \hat{h}^k(\tilde{x}^{k+2})$, we
get \BEAS 
\delta^k - \delta^{k+1} 
&\geq&\hat{h}^{k+1}(\tilde{x}^{k+2})-\hat{h}^k\left(\tilde{x}^{k+1}
\right) - \rho\left( x^p - \tilde{x}^{k+1}\right)^T(\tilde{x}^{k+2}
-\tilde{x}^{k+1})\\
&& + (\rho/2)\left\|\tilde{x}^{k+2} -\tilde{x}^{k+1}\right \|_2^2\\
&\geq& (\rho/2)\left\|\tilde{x}^{k+2} -\tilde{x}^{k+1}\right \|_2^2. 
\EEAS
Therefore, $\delta^k \geq \delta^{k+1}+ (\rho/2)\left\|\tilde{x}^{k+2} -
\tilde{x}^{k+1}\right \|_2^2$ for all $k \geq p+1$. Then from 
\BEAS
\hat{h}^k(x^p)&\geq&  \hat{h}^k\left(\tilde{x}^{k+1}\right) +
\rho(x^p - \tilde{x}^{k+1})^T(x^p - \tilde{x}^{k+1})\\
&=&h(x^k)-\delta^k+(\rho/2)\left \|x^p - \tilde{x}^{k+1}\right
\|_2^2,
\EEAS
it follows that $\left\|x^p - \tilde{x}^{k+1}\right \|_2^2\leq
2\delta^k/\rho\leq 2\delta^p/\rho$.

Now we use the assumption that
$f$ and $g$ are Lipschitz continuous with Lipschitz constant $L$ for all $ x
\in \dom g$. Every
$q\in \partial \hat{f}^k(x)$ has the form $q= \sum_{t \leq k} \theta_t q^t$,
with $\theta_t \geq 0$ and $\sum_t \theta_t=1$,  a convex
combination of normal vectors of active constraints at $x$, where $q^t \in
\partial f(x^t)$. Therefore, $\hat{h}^k(x)=\hat{f}^k (x)+g(x)$
is $2L$-Lipschitz continuous.

Combining this with 
\[
 \delta^k \leq h(x^k) -
 \hat{h}^k\left(\tilde{x}^{k+1}\right), \qquad -\eta\delta^k \leq
 {h}\left(\tilde{x}^{k+1}\right) - h(x^k), 
\]  
we have 
\[
 (1-\eta)\delta^k \leq {h}\left(\tilde{x}^{k+1}\right)  -
 h\left(\tilde{x}^k\right) + \hat{h}^k\left(\tilde{x}^k\right)  -
 \hat{h}^k\left(\tilde{x}^{k+1}\right)
\leq 4L \left\|\tilde{x}^k  -  \tilde{x}^{k+1} \right\|_2.
\]
Therefore, from 
\BEAS 
\frac{(1-\eta)^2\rho}{32L^2} \sum_{k\geq p}
\left(\delta ^k\right)^2 \leq\sum_{k\geq p}
\left(\delta^k-\delta^{k+1}\right) \leq \delta^{p}, 
\EEAS
we can establish that $\delta^{k}$ converges to zero as $k \to \infty$. 
This implies
\[
\underset{k\to \infty}{\lim}\left(\hat{h}^k\left(\tilde{x}^{k+1}\right)
+(\rho/2) \left\| \tilde{x}^{k+1} - x^p\right\|_2^2 \right)= h(x^p).
\]
Also from $\underset{k\to \infty}{\lim} \left (h(x^p) -
h\left(\tilde{x}^{k+1}\right) \right) = 0$ and  $\left\|{x}^{p} -
\tilde{x}^{k+1}\right \|_2^2 \leq \frac{2\delta^k}{\rho}$, it follows that 
\[
\underset{k\to \infty}{\lim} \hat{h}^k\left(\tilde{x}^{k+1}
\right) = h(x^p),
\qquad\underset{k\to \infty}{\lim} \left\|
{x}^{p} - \tilde{x}^{k+1}\right\|_2^2 = 0.
\]
Hence, we get $0 \in \partial h(x^p)$, which implies $h(x^p)=h^\star$.

\end{document}